\documentclass{elsart}  \usepackage{amsfonts} \usepackage{amssymb} \usepackage{latexsym} \usepackage{graphicx} \usepackage{hyperref} \usepackage{amsbsy} \usepackage{psfrag}       \newcommand{\D}{\displaystyle}
\newcommand{\restr}[1]{\raisebox{-0.3em}{$\lb|_{#1}\rb.$}} 
\newcommand{\ignore}[1]{}    
\newcommand{\breath}{\medskip} 
\newtheorem{THM}{Theorem}[section]

\newcounter{claimcount}[THM]  
\newcounter{subclaimcount}[claimcount]
\newtheorem{PROP}[THM]{Proposition} 

\newtheorem{lemma}[THM]{Lemma} 

\newtheorem{COR}[THM]{Corollary}

\newcommand{\dfn}{\sf\em} 
\newcommand{\Theorem}[2]{\begin{THM}{\sf #1}  #2 \end{THM}}
\newcommand{\Proposition}[2]{\begin{PROP}{\sf #1}  #2 \end{PROP}}
\newcommand{\Lemma}[2]{\begin{lemma}{\sf #1}  #2 \end{lemma}}
\newcommand{\Corollary}[2]{\begin{COR}{\sf #1}  #2 \end{COR}} 
\newcommand{\THMfont}[1]{{\sl #1}}      
\newcommand{\bthmlist}{ \begin{list}{{\bf(\alph{enumi})}} {\usecounter{enumi} \setlength{\itemsep}{0.2em} \setlength{\topsep}{0.2em} \setlength{\itemindent}{0em} \setlength{\parsep}{0em} } } 
\newcommand{\ethmlist}{\end{list}}    
\newcommand{\Claim}[1]{\refstepcounter{claimcount} \vspace{0.3em}               \noindent {\sc Claim \theclaimcount: \ }\THMfont{ #1}} 
\newcommand{\Subclaim}[1]{\refstepcounter{subclaimcount} \vspace{0.3em}               \noindent {\sc Claim \theclaimcount.\thesubclaimcount: \ }\THMfont{ #1}} 
\newcommand{\subclaim}{\Subclaim}
\newcommand{\bprf}[1][Proof:]{\begin{list}{} 			{\setlength{\leftmargin}{0.5em} 			\setlength{\rightmargin}{0em} 			\setlength{\listparindent}{1em}}                         \item {\em \hspace{-0.8em}  #1  }} 
\newcommand{\eprf}{\end{list}} 
\newcommand{\bthmprf}{\bprf}
\newcommand{\bclaimprf}{\bprf}
\newcommand{\bsubclaimprf}{\bprf}
\newcommand{\eclaimprf}{ \hfill $\Diamond$~{\scriptsize {\tt Claim~\theclaimcount}}\eprf}  
\newcommand{\esubclaimprf}{ \hfill $\triangledown$~{\scriptsize  {\tt Claim~\theclaimcount.\thesubclaimcount}}\eprf}  

\newcommand{\beq}{\begin{eqnarray*}}
\newcommand{\eeq}{\end{eqnarray*}} 
\newcommand{\beqn}{ \begin{equation} }
\newcommand{\eeqn}{ \end{equation} }
\newcommand{\bitem}{\begin{itemize}}
\newcommand{\eitem}{\end{itemize}} 
\newcommand{\bdesc}{\begin{description}}
\newcommand{\edesc}{\end{description}}   
\newcommand{\If}{\mbox{\ if \ }}  
\newcommand{\done}{{\mathsf{ 1\!\!1}}} 
\newcommand{\dA}{{\mathbb{A}}}
\newcommand{\dB}{{\mathbb{B}}}
\newcommand{\dE}{{\mathbb{E}}}
\newcommand{\dH}{{\mathbb{H}}}
\newcommand{\dN}{{\mathbb{N}}}
\newcommand{\dR}{{\mathbb{R}}}
\newcommand{\dZ}{{\mathbb{Z}}}       
\newcommand{\barB}{{\overline{B}}}
\newcommand{\barS}{{\overline{S}}}
\newcommand{\baralp }{{\overline{\alpha}}}
\newcommand{\barba}{{\overline{\mathbf{ a}}}}
\newcommand{\bA}{{\mathbf{ A}}}
\newcommand{\bB}{{\mathbf{ B}}}
\newcommand{\bC}{{\mathbf{ C}}}
\newcommand{\bD}{{\mathbf{ D}}}
\newcommand{\bH}{{\mathbf{ H}}}
\newcommand{\bI}{{\mathbf{ I}}}
\newcommand{\bK}{{\mathbf{ K}}}
\newcommand{\bL}{{\mathbf{ L}}}
\newcommand{\bO}{{\mathbf{ O}}}
\newcommand{\bS}{{\mathbf{ S}}}
\newcommand{\bT}{{\mathbf{ T}}}
\newcommand{\bU}{{\mathbf{ U}}}
\newcommand{\bV}{{\mathbf{ V}}}
\newcommand{\bW}{{\mathbf{ W}}}
\newcommand{\bX}{{\mathbf{ X}}}
\newcommand{\bY}{{\mathbf{ Y}}}
\newcommand{\ba}{{\mathbf{ a}}}
\newcommand{\bb}{{\mathbf{ b}}}
\newcommand{\bc}{{\mathbf{ c}}}
\newcommand{\bk}{{\mathbf{ k}}}
\newcommand{\bo}{{\mathbf{ o}}}
\newcommand{\br}{{\mathbf{ r}}}
\newcommand{\bu}{{\mathbf{ u}}}
\newcommand{\sA}{{\mathcal{ A}}}
\newcommand{\sB}{{\mathcal{ B}}}
\newcommand{\sC}{{\mathcal{ C}}}
\newcommand{\sE}{{\mathcal{ E}}}
\newcommand{\sF}{{\mathcal{ F}}}
\newcommand{\sK}{{\mathcal{ K}}}
\newcommand{\sM}{{\mathcal{ M}}}
\newcommand{\sO}{{\mathcal{ O}}}
\newcommand{\sR}{{\mathcal{ R}}}
\newcommand{\gA}{{\mathfrak{ A}}}
\newcommand{\gG}{{\mathfrak{ G}}}
\newcommand{\gK}{{\mathfrak{ K}}}
\newcommand{\gL}{{\mathfrak{ L}}}
\newcommand{\gO}{{\mathfrak{ O}}}
\newcommand{\gb}{{\mathfrak{ b}}}
\newcommand{\gs}{{\mathfrak{ s}}}
\newcommand{\alp }{\alpha}
\newcommand{\bet }{\beta}
\newcommand{\gam }{\gamma}
\newcommand{\del }{\delta}
\newcommand{\eps }{\epsilon}
\newcommand{\kap }{\kappa}
\newcommand{\lam }{\lambda}
\newcommand{\sig }{\sigma} 
\newcommand{\ups }{\upsilon}
\newcommand{\Del }{\Delta}
\newcommand{\fk}{{\mathsf{ k}}}
\newcommand{\fv}{{\mathsf{ v}}}
\newcommand{\fw}{{\mathsf{ w}}}
\newcommand{\fz}{{\mathsf{ z}}}    
\newcommand{\tla}{{\widetilde{a}}}
\newcommand{\tlba}{{\widetilde{\mathbf{ a}}}}
\newcommand{\undB}{{\underline{B}}}
\newcommand{\undS}{{\underline{S}}}
\newcommand{\vV}{{\vec{V}}}
\newcommand{\lb}{\left}
\newcommand{\rb}{\right} 
\newcommand{\maketall}{\rule[-0.5em]{0em}{1em}}        
\newcommand{\map}{{\longrightarrow}}
\newcommand{\goto}{{\rightarrow}}
\newcommand{\into}{{\map}}
\newcommand{\statement}[1]{\lb(  \maketall       \begin{minipage}{40em}       \begin{tabbing}         #1        \end{tabbing}      \end{minipage}  \rb)}     
\newcommand{\oo}{{\infty}}        
\newcommand{\X}{\times}
\newcommand{\x}{\X}
\newcommand{\compl}[1]{#1^\complement}  
\newcommand{\symdif}{{\bigtriangleup}} 
\newcommand{\union}{\cup}
\newcommand{\Union}{\bigcup}
\newcommand{\intsct}{\cap}
\newcommand{\Intsct}{\bigcap}
\newcommand{\disj}{\sqcup}
\newcommand{\Disj}{\bigsqcup}   
\newcommand{\set}[2]{{\left\{ #1 \; ; \; #2 \right\} }} 
\newcommand{\supp}[1]{{\sf supp}\lb(#1\rb)}      
\newcommand{\norm}[2]{{\left\| #1 \right\|_{{#2}} }   }
\newcommand{\chr}[1]{{{\done}_{{#1}}}} 
\newcommand{\choice}[1]{{\lb\{ \begin{array}{rcl}                                 #1                                \end{array}  \rb.  }}                     %
\newcommand{\eeequals}[1]{\raisebox{-0.78ex}{$\overline{\overline{{\scriptscriptstyle{\mathrm{#1}}}}}$}} 
\newcommand{\leeeq}[1]{\raisebox{-1ex}{${{\D\leq} \atop {\scriptscriptstyle{\mathrm{#1}}}}$}} 
\newcommand{\lt}[1]{\raisebox{-1ex}{${{\D<} \atop {\scriptscriptstyle{\mathrm{#1}}}}$}} 
\newcommand{\suuubset}[1]{\raisebox{-1.3ex}{$\stackrel{\D\subset}{\scriptscriptstyle{\mathrm{#1}}}$}}  
\newcommand{\iiimplies}[1]{\eeequals{#1}\!\!\!\!\!\!\Rightarrow}
\newcommand{\closeto}[1]{{{\raisebox{-1ex}    {$\widetilde{\ {\scriptstyle #1}\ }$}}}}      
\newcommand{\cl}[1]{{\sf cl}\lb(#1\rb)}
\newcommand{\shift}[1]{\sig^{#1}}    
\newcommand{\goesto}[2]{{ -\!\!\!-\!\!\!-\!\!\!-\!\!\!\!\!\!\!\!\!\!\!  ^{{\scriptscriptstyle #2}}_{{\scriptscriptstyle #1}}   \!\!\!\!\!\!\!\!\!\longrightarrow }}                         
\newcommand{\diam}[1]{{\sf diam}\lb[#1\rb]} 
\newcommand{\grad}{{\nabla}}
\newcommand{\Real}{\dR}
\newcommand{\Natur}{\dN}
\newcommand{\Zahl}{\dZ}
\newcommand{\CC}[1]{{\lb[ #1 \rb]}}
\newcommand{\CO}[1]{{\lb[ #1 \rb)}}
\newcommand{\OC}[1]{{\lb( #1 \rb]}}
\newcommand{\OO}[1]{{\lb( #1 \rb)}}   
\newcommand{\dX}{{\Real^D}}
\newcommand{\dK}{\bK}
\newcommand{\gk}{\kappa}  
\newcommand{\eref}[1]{eqn.{\rm(\ref{#1})}}           
\newcommand{\Ll}{\bL^1}
\newcommand{\Loo}{\bL^\oo}
\newcommand{\elll}{\ell^1}  
\newcommand{\Lllim}{\Ll\!\!\!-\!\!\!\lim}
\newcommand{\Loolim}{\Loo\!\!-\!\!\lim}        
\newcommand{\Life}{\Upsilon}      
\newcommand{\Int}[1]{\mathsf{int}\lb(#1\rb)} 
\newcommand{\HD}[1]{d_{\mathrm{H}}\lb(#1\rb)}  
\newcommand{\tlgk}{\tilde{\gk}}
\newcommand{\bargk}{\bar{\gk}}  
\newcommand{\AX}{\sA^\dX}
\newcommand{\lAX}{{}^1\!\!\AX}
\newcommand{\bAX}{{}^\partial\!\!\AX}
\newcommand{\oAX}{{}^0\!\!\AX}
\newcommand{\dAX}{{}^*\!\AX}  
\newcommand{\AZD}{\sA^{\Zahl^D}}
\newcommand{\lAZD}{{}^1\!\!\AZD}
\newcommand{\ZD}{\Zahl^D}
\newcommand{\eZD}[1][\eps]{#1\Zahl^D}          
\newcommand{\eba}[1][\eps]{\barba_{#1}}
\newcommand{\ebgk}[1][\eps]{\bargk_{#1}}
\newcommand{\etgk}[1][\eps]{\tlgk_{#1}}
\newcommand{\eLl}[1][\eps]{\Ll_{#1}}
\newcommand{\eAX}[1][\eps]{{}^{#1}\!\AX}
\newcommand{\eLife}[1][\eps]{\Life_{#1}}
\newcommand{\etLife}[1][\eps]{\widetilde{\Life}_{#1}}
\newcommand{\ebLife}[1][\eps]{\bar{\Life}_{#1}}    
\newcommand{\ethmprf}{\hfill$\Box$\eprf}
\newcommand{\implies}{\ensuremath{\Longrightarrow}}\newcommand{\And}{\mbox{\ and \ }}  

\begin{document}

\begin{frontmatter}
\title{RealLife:  the continuum limit of Larger Than Life cellular automata}
\author{Marcus Pivato}
\address{Department of Mathematics,  Trent University\\
 1600 West Bank Drive, Peterborough, Ontario, Canada, K9J 7B8}

\ead{pivato@xaravve.trentu.ca}

\begin{abstract}
  Let $\sA:=\{0,1\}$.  A {\em cellular automaton} (CA) is a
shift-commuting transformation of $\sA^{\Zahl^D}$ determined by a
local rule. Likewise, a {\em Euclidean automaton} (EA) is a
shift-commuting transformation of $\sA^{\Real^D}$ determined by a
local rule.  {\em Larger than Life} (LtL) CA are long-range
generalizations of J.H. Conway's {\em Game of Life} CA,
proposed by K.M. Evans.  We prove a conjecture of Evans: as their radius
grows to infinity, LtL CA converge to a `continuum limit' EA, which we
call {\em RealLife}.  We also show that the {\em life forms} (fixed
points, periodic orbits, and propagating structures) of LtL CA
converge to life forms of {\em RealLife}.
Finally we prove a number of existence results for fixed points of
{\em RealLife}.

{\breath

\begin{tabular}{rl}
{\bf MSC:}& 37B15 (primary), 68Q80 (secondary)\\
{\bf Keywords:}& Cellular automata, Game of Life, continuum limit.
\end{tabular}}

\end{abstract}

\end{frontmatter}


  Let $\sA:=\{0,1\}$ and let $D\in\Natur$.  Let $\AZD$ be the set
of all {\dfn configurations} ---ie. functions $\ba:\ZD\into\sA$.
If $\dK\subset\ZD$, then we define $\ba\restr{\dK}\in\sA^\dK$
to be the restriction of $\ba$ to $\dK$.
A $D$-dimensional {\dfn cellular automaton} (CA) is a transformation 
$\Phi:\AZD\into\AZD$ determined by a finite subset $\dK\subset\ZD$
 (the {\dfn neighbourhood}) and a {\dfn local rule} $\phi:\sA^\dK\into\sA$
so that, for any $\ba\in\AZD$, $\Phi(\ba):=\bb$, where
$\bb(\fz) \ = \ \phi(\ba\restr{\fz+\dK})$ for all $\fz\in\ZD$.
 
  One of the most  fascinating cellular automata is
Conway's {\dfn Game of Life} \cite{BerlekampConwayGuy}, which  is the 
function
$\Life:\sA^{\Zahl^2}\into\sA^{\Zahl^2}$ defined
\[
  \Life(\ba)(\fz) \  := \
\choice{1 &&\If  \ba(\fz) = 1 \And  \gK(\ba)(\fz)\in\{3,4\}; \\
        1 &&\If  \ba(\fz) = 0 \And  \gK(\ba)(\fz) = 3; \\
        0 &&\mbox{otherwise.}}
\]
Here, $\D\gK(\ba)(\fz) \ := \ \D \sum_{\fk\in\dK} \ba(\fz+\fk)$, \ 
where $\dK  \ := \ \CC{-1...1}\x\CC{-1...1}$.

\begin{figure}
\begin{center}
\begin{tabular}{cccc}
\includegraphics[scale=0.8]{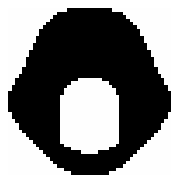} &
\includegraphics[scale=0.8]{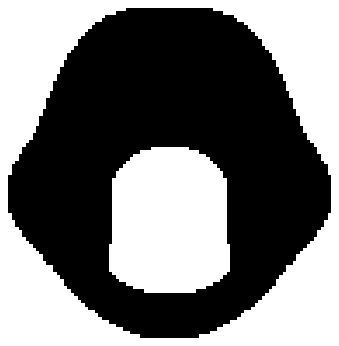}&
\includegraphics[scale=0.8]{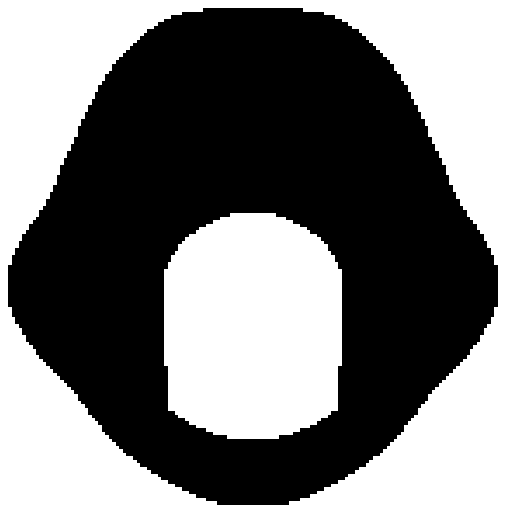}&
\includegraphics[scale=0.8]{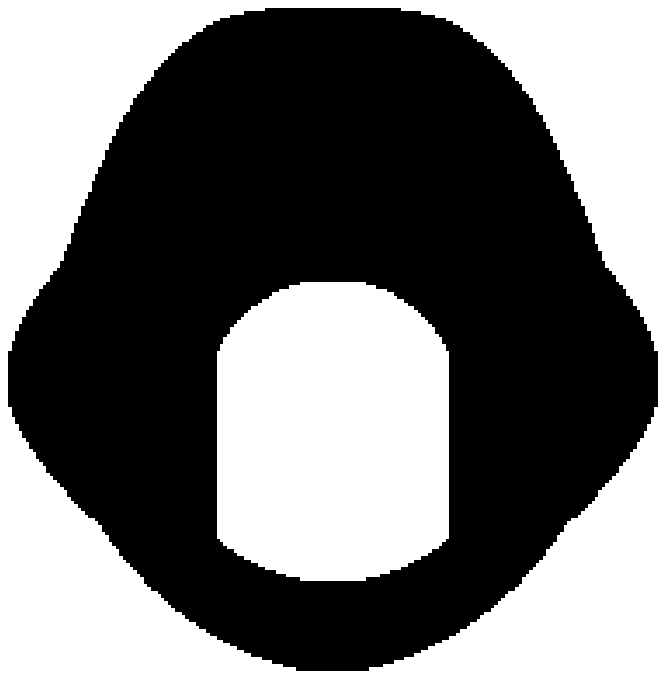}\\
K=25 & K=50 & K=75 & K=100\\ 
\end{tabular}
\end{center}
\caption{{\footnotesize Four morphologically similar bugs for LtL CA of increasing
radius. In all cases, $s_0=b_0 = \frac{706}{2601} < b_1 = \frac{958}{2601} < s_1 = \frac{1216}{2601}$,
and $\dK=\CC{-K...K}^2$, where $K=25$, $50$, $75$, or $100$.
\label{fig:four.bugs}}}
\end{figure}

  {\dfn Larger than Life} (LtL) is an infinite family of long-range,
two-dimensional generalizations of {\em Life}, introduced
by Evans \cite{Evans1,Evans2,Evans3,Evans4,Evans5}.  An LtL CA has the
form:
\[
  \Life(\ba)(\fz) \ := \
\choice{1 &&\mbox{if} \  \ba(\fz) = 1 \And  s_0\leq \gK(\ba)(\fz) \leq s_1; \\
        1 &&\mbox{if} \  \ba(\fz) = 0 \And  b_0\leq \gK(\ba)(\fz) \leq b_1; \\
        0 &&\mbox{otherwise.}}
\]
  Here  $0 \leq s_0 \leq b_0 \leq b_1 \leq s_1 \leq 1$, and
$\gK(\ba)(\fz) \ := \ \D \frac{1}{|\dK|}\sum_{\fk\in\dK} \ba(\fz+\fk)$, where
$\dK$ is some large `neighbourhood' of the origin, usually
$\dK=\CC{-K...K}\x\CC{-K...K}$ (for some radius $K>0$).
 More generally, we could define
$\gK(\ba)(\fz) := \sum_{\fk\in\dK} c_\fk \ba(\fz+\fk)$ for any set of
nonnegative coefficients $\{c_{\fk}\}_{\fk\in\dK}$ 
such that $\sum_{\fk\in\dK} c_\fk=1$.
 We refer to $\CC{b_0,b_1}$ as the
{\dfn birth interval} and $\CC{s_0,s_1}$ as the
{\dfn survival interval}.  In Conway's {\em Life},
$s_0=b_0=b_1=\frac{1}{3}$, and $s_1=\frac{4}{9}$.  
Evans' {\em Larger than Life} CA usually have
\[
0.2  \ \leq \  s_0 \ \leq \ b_0 \ \leq \  0.27
 \ \leq \  0.3  \ \leq \ b_1 \ \leq \ 0.35 \  \leq \ s_1 \  \leq \ 0.5.
\] 
 {\em Larger than Life} CA exhibit phenomena qualitatively
similar to those found in {\em Life} and its generalizations 
(see \cite{Bays6,Bays5,Bays4,Bays3,Bays2,Bays1}, \cite{Eppstein,Gotts}
and \cite[\S6]{GrGr})
including the emergence of
complex, compactly supported fixed points ({\dfn still lifes}),
periodic solutions ({\dfn oscillators}) and propagating structures
called {\dfn bugs} (analogous to the {\em gliders} and {\em
space-ships} of {\em Life}), which can sometimes be arranged to 
perform computation \cite{Evans5}.  Especially
intriguing is that the still lifes, oscillators, and bugs found in
longer-range LtL CA appear to be rescaled, `high resolution'
versions of those found in shorter range LtL CA (see Figure~\ref{fig:four.bugs}).  Evans \cite{Evans2} conjectures that these
still lifes (resp. oscillators, bugs) converge to a continuum limit,
which is a still life (resp. oscillator, bug) for
some kind of {\dfn Euclidean automaton}; a
translationally-equivariant transformation of $\sA^{\Real^2}$.

\ignore{{\bf(A)} An eastward moving bug in the range 5 LtL CA 
with $s_0=b_0 = \frac{24}{121} < b_1 = \frac{30}{121} < s_1 = \frac{33}{121}$.
\quad
{\bf(B)} A highly symmetric still life in the range 5 LtL CA 
with $s_0=b_0 = \frac{25}{121} < b_1 = \frac{32}{121} < s_1 = \frac{35}{121}$.}

\breath

  In \S\ref{S:euclid}, we formally define Euclidean automata (EA), and
introduce the {\em RealLife} family of EA, the natural generalization
of {\em Larger than Life} to $\sA^{\Real^D}$.  We show that {\em
RealLife} EA are continuous on a comeager set in the natural
$\bL^1$ norm on $\sA^{\Real^D}$ (Theorems
\ref{reallife.is.continuous} and \ref{oAX.is.comeager}), and show that their dynamics vary continuously
as a function of the parameters $(s_0,b_0,b_1,s_1)$ and the neighbourhood
$\dK$ (Theorems \ref{threshold.convergence} and \ref{akylea.convergence}).  In
\S\ref{S:LtL} we show that {\em RealLife} EA are the
continuum limits of suitable sequences of {\em Larger than Life} CA of
increasing radius, and that a suitable converging
sequence of still lifes (resp. oscillators, bugs) for these LtL CA
yields a still life (resp. oscillator, bug) for {\em RealLife}
(Theorem \ref{pointwise.convergence}).  In
\S\ref{S:rot.fixed.point}, we construct several families of nontrivial
still lifes for {\em RealLife} EA satisfying various conditions.
Finally, in \S\ref{S:hausdorff}, we introduce the {\em Hausdorff metric} $d_*$
on $\sA^{\Real^D}$, and show that a still life is often surrounded by
a $d_*$-neighbourhood of other still lifes (Theorem
\ref{fixed.point.nhood}).

The four sections are mostly logically independent, except for the
use of notation and definitions from \S\ref{S:euclid}.
Also, the proof of Theorem \ref{pointwise.convergence} uses
Theorem \ref{reallife.is.continuous}, Lemma \ref{conv.lemma.2} and
Proposition \ref{fixed.point.convergence.lemma}, and the proof of 
Theorem \ref{fixed.point.nhood} uses Proposition \ref{fixed.point.cond}.

\section{Euclidean Automata and RealLife
\label{S:euclid}}

Let $\lam$ be the $D$-dimensional Lebesgue
measure on $\dX$, and let $\Loo := \Loo(\dX,\lam)$.  Let
$\sA:=\{0,1\}$ and let $\AX\subset\Loo$ be the set of all
Borel-measurable functions $\ba:\dX\into\sA$, which we will refer to
as  {\dfn configurations}.  If $v\in\dX$, then define the shift map
$\shift{v}:\AX\into\AX$ by $\shift{v}(\ba)=\ba'$, where $\ba'(x) =
\ba(x+v)$ for all $x\in\dX$.  A {\dfn Euclidean automaton} (EA) is a
function $\Phi:\AX\into\AX$ which commutes with all shifts, and which
is {\em determined by local information}, meaning that there is some
compact neighbourhood $\dK\subset\dX$ around zero so that, if
$\ba,\ba'\in\AX$, and $\ba\restr{\dK}=\ba'\restr{\dK}$, then
$\Phi(\ba)(0)=\Phi(\ba')(0)$.

  Let $\sK:=\set{\gk\in\bL^\oo(\dX;\CO{0,\oo})}{ \gk \ \mbox{has compact support,
 and} \ \int_\dX \gk \ = \ 1}$.  If $\gk\in\sK$ and 
$\ba\in\sA^\dX$, then the convolution of $\ba$ by $\gk$ is defined:
\[
 \gk*\ba(x) \quad := \quad \int_{\dX} \gk(y) \cdot \ba(x-y) \ d\lam[y].
\]
For example, 
if $\dK\subset\dX$ is a compact neighbourhood of zero (eg. a ball or a cube),
 and $\gk:=\lam[\dK]^{-1} \chr{\dK}$, then $\gk\in\sK$, and
 $\gk*\ba(x)  \ = \ \lam[\dK]^{-1} \int_{\dK} \ba(x-k) \ d\lam[k]$
 is the average value of $\ba$ near $x$.
If $0 < s_0 \leq b_0 < b_1 \leq s_1\leq 1$, then
the corresponding {\dfn RealLife} Euclidean automaton
$\Life={}_\gk\Life_{b_0,s_0}^{b_1,s_1}:\sA^\dX\into\sA^\dX$ is defined:
\beqn
\label{afterlife.CA1}
\hspace{-2em}
\forall \ \ba\in\AX, \quad
  \Life(\ba)(x) \quad := \quad
\choice{1 &&\If  \ba(x) = 1 \And \ s_0\leq \gk*\ba(x) \leq s_1; \\
        1 &&\If  \ba(x) = 0 \And \ b_0\leq \gk*\ba(x) \leq b_1; \\
        0 &&\mbox{otherwise.}}\qquad
\eeqn
Let $\Theta:=\set{(s_0,b_0,b_1,s_1)}{0< s_0 \leq b_0 < b_1 \leq s_1\leq 1}$
be the set of {\dfn threshold four-tuples}.
Note that $\Life$ depends upon the choice of kernel $\gk\in\sK$ and the
four-tuple $(s_0,b_0,b_1,s_1)\in\Theta$;  we normally 
suppress this dependency in our notation.
We call this Euclidean automaton {\em RealLife} because it is
the continuum limit of a sequence of
{\em Larger than Life} cellular automata
with  `birth interval' $\CC{b_0,b_1}$ and `survival interval'
$\CC{s_0,s_1}$ (Theorem \ref{pointwise.convergence}).
If we define $\gb:=\chr{\CC{b_0,b_1}}$ and $\gs := \chr{\CC{s_0,s_1}}$,
then, for any $\ba\in\AX$ and $x\in\dX$,
we can rewrite eqn.(\ref{afterlife.CA1}) as
\beqn
\label{afterlife.CA2}
  \Life(\ba)(x) \quad := \quad  \ba(x)\cdot \gs\lb(\gk*\ba(x)\rb) 
 \ + \  (1-\ba(x))\cdot \gb\lb(\gk*\ba(x)\rb).
\eeqn
  The {\dfn compact-open} topology on $\AX$
is determined by the metric $d_C$ defined for all $\ba,\ba'\in\AX$ by
\[
\hspace{-2em}
  d_C(\ba,\ba')\ := \ \exp\lb[-R(\ba,\ba')\rb], \ 
\mbox{where} \  R(\ba,\ba') \ := \ 
\sup\set{r>0}{\ba\restr{\dB(r)} \equiv  \ba'\restr{\dB(r)}}.
\]
 (Here $\dB(r):=\set{x\in\dX}{|x|\leq r}$).  It is not hard to prove the 
analog of the Curtis-Hedlund-Lyndon Theorem 
\cite{HedlundCA} that is  well-known for cellular automata:

\paragraph*{Theorem:}
{\em
  Let $\Phi\colon\AX\goto\AX$ be $\shift{}$-commuting.  Then $\Phi$ is
an EA iff $\Phi$ is $d_C$-continuous.\qed
}

\breath

  The compact-open topology is very fine; $d_C(\ba,\ba')$ will be
large even if $\ba$ and $\ba'$ differ on a set $\Del\subset\dX$ of
tiny measure, as long as $\Del$ contains points near the origin.
Also, $d_C$ is not shift-invariant.  Thus, sometimes it is
more suitable to use the $\bL^1$ metric.  Let $\Ll:=\Ll(\dX,\lam)$ and
let $\lAX := \AX\intsct\Ll$ be the set of configurations whose support
has finite measure. 
Note that $\Life(\lAX)\subseteq\lAX$ (because $0<b_0$).
We extend $\Life$ to a function $\Life:\Ll\into\Ll$
by applying eqn.(\ref{afterlife.CA2}) in the obvious way.
For any $\ba\in\lAX$, 
let $M(\ba):= \lam\lb[\alp^{-1}\{s_0,b_0,s_1,b_1\}\rb]$, where
$\alp:=\gk*\ba$. We define
\[
  \oAX\quad:=\quad\set{\ba\in\lAX}{M(\ba)=0}.
\]
(Note that $\oAX$ is a function of $(s_0,b_0,s_1,b_1)$ and
$\gk$).  This section's  first main result  is:

\Theorem{\label{reallife.is.continuous}}
{
If $(s_0,b_0,b_1,s_1)\in\Theta$ and $\gk\in\sK$, then
 $\Life$ is $\Ll$-continuous on $\oAX$.
}

 $\oAX$ is a strict subset of $\lAX$.  To see this, suppose
$\gk=\chr{\dK}$, where
$\dK=\CC{\frac{-1}{2},\frac{1}{2}}^D$.  Let $r:=\sqrt[D]{s_0}$ and let
$\bA := \CC{0,r}^D$, so $\lam[\bA]=s_0$.  If $\ba:=\chr{\bA}$, and
 $\alp:=\gk*\ba$, then $\alp(x)=s_0$ for $\forall \ x\in\CC{r-1,1}^D$.
Thus, $M(\ba)=\lam\lb(\CC{r-1,1}^D\rb) = (2-r)^D>0$, so $\ba\not\in\oAX$.

 $\Life$ is {\em not} $\bL^1$-continuous on all of $\lAX$.  To see this,
 note that $\Life(\ba)=\chr{\bB}$, where
$\bB:=\bA$ if $b_0>s_0$, and $\bB :=\CC{r-1,1}^D\supset \bA$ if
$b_0=s_0$. Now, let $\eps>0$ be tiny. Let $r':=r-\eps/D$, let $\bA' :=
\CC{0,r'}^D$, and let $\ba':=\chr{\bA'}$. Then $\norm{\ba-\ba'}{1}<\eps$.
However,  $\lam[\bA']<s_0$, so
that $\Life(\ba')=\bo$ is the zero configuration.  Thus,
$\norm{\Life(\ba)-\Life(\ba')}{1} = \lam[\bB] \geq s_0$.

  Fortunately, this discontinuity set is usually meager in
$\lAX$.  If $\bT\subset\dX$, and $\gam>0$,
let $\dB(\bT,\gam):=\set{x\in\dX}{d(x,\bT)<\gam}$.
Say $\bT$ is {\dfn thin} if $\D\lim_{\gam\goto0} \lam[\dB(\bT,\gam)]=0$.
For example, any compact, piecewise smooth $(D-1)$-submanifold of $\dX$ is
thin.  The kernel $\gk$ is {\dfn almost continuous}
if there is a thin set $\bT\subset\dX$ so that, for any $\gam>0$,
$\gk$ is uniformly continuous on $\compl{\dB(\bT,\gam)}$.
For example:
\bitem
  \item If $\gk$ is
continuous, then $\gk$ is almost continuous ($\gk$ has compact support, so
continuity implies uniform continuity).
  
  \item If $\dK$ is an open
set and $\partial\dK$ is thin (eg. $\partial\dK$ is a piecewise
smooth manifold), then $\gk:=\lam[\dK]^{-1}\chr{\dK}$ is almost continuous.
\eitem

\Theorem{\label{oAX.is.comeager}}
{
For any $(s_0,b_0,b_1,s_1)\in\Theta$ and any
almost-continuous $\gk\in\sK$, the set \
$\oAX$ is a $\shift{}$-invariant dense G$\delta$ subset of $\lAX$.
}

Let $\Phi:\Ll\into\Ll$ be a $\shift{}$-commuting transformation. If
$\ba\in\Ll$, then $\ba$ is a {\dfn still life} for $\Phi$ if
$\Phi(\ba)=\ba$.  If $p\in\Natur$, then $\ba$ is a {\dfn
$p$-oscillator} if $\Phi^p(\ba)=\ba$ (a still life is thus a
$1$-oscillator).  If $p\in\Natur$ and $v\in\dX$, then $\ba$ is a {\dfn
$p$-periodic bug} with {\dfn velocity $v$} if $\Phi^p(\ba) \ = \
\shift{pv}(\ba)$.  We will refer to still lifes, oscillators, and bugs
collectively as {\dfn life forms}.

  Recall that $\Life$ is determined by the threshold parameter four-tuple
$(s_0,b_0,b_1,s_1)$ and the convolution kernel $\gk$.  A small change
in these parameters should yield a small change in $\Life$, and a
small change in its life forms.  In particular, if
$\gk$ is a fixed kernel, and
$\{(s^n_0,b^n_0,b^n_1,s^n_1)\}_{n=1}^\oo$ is a sequence of 
four-tuples converging to the four-tuple $(s_0,b_0,b_1,s_1)$, then the
corresponding sequence $\{\Life_n\}_{n=1}^\oo$ of RealLife EA (with
kernel $\gk$) should converge to the RealLife EA $\Life$
determined by $(s_0,b_0,b_1,s_1)$ and $\gk$.  Likewise, if we fix
$(s_0,b_0,b_1,s_1)$, then a convergent sequence 
$\{\gk_n\}_{n=1}^\oo\subset\sK$ of  kernels
should yield a convergent sequence of RealLife EA.  Furthermore, in both
cases, the life forms of $\{\Life_n\}_{n=1}^\oo$ should
 `evolve' toward life forms for $\Life$.

  To formalize life form `evolution', suppose
$\{\Phi_n\}_{n=1}^\oo$ was a sequence of $\shift{}$-commuting
transformations of $\Ll$.  If $\gA\subset\Ll$ is a
$\shift{}$-invariant subset, then the sequence $\{\Phi_n\}_{n=1}^\oo$
{\dfn evolves} to $\Phi$ on $\gA$ if, for any
$\{\ba_n\}_{n=1}^\oo\subset\Ll$ such that $\D\Lllim_{n\goto\oo} \ba_n
\ = \ \ba\in\gA$, the following holds:
 \bthmlist
  \item If $\Phi_n(\ba_n)=\ba_n$ for all $n\in\Natur$, 
 then $\Phi(\ba) \ = \ \ba$.

  \item Let $P\in\Natur$, and suppose
 $\Phi^p(\ba)\in\gA$ for all $p\in\CO{0...P}$.  
\bitem
  \item[{\bf[i]}] If $\Phi_n^P(\ba_n)=\ba_n$ for all $n\in\Natur$, 
 then $\Phi^P(\ba) \ = \ \ba$.

 \item[{\bf[ii]}] If $\{v_n\}_{n=1}^\oo\subset\dX$
and $\D\lim_{n\goto\oo} v_n = v\in\dX$, and
 $\Phi_n^P(\ba_n)=\shift{Pv_n}(\ba_n)$ for all $n\in\Natur$, 
 then $\Phi^P(\ba) \ = \ \shift{Pv}(\ba)$.
\eitem
\ethmlist
\ignore{  To get `evolution' for a convergent sequence
of threshold four-tuples, we need another technicality.
  If $\ba\in\lAX$, and $\alp=\gk*\ba$, then
$\mu_\ba := \alp^*(\lam)$ is a measure on $\CC{0,1}$. 
Let $\lam_1$ be the Lebesgue measure on $\Real$, and
define
\[
  M'(\ba) \ := \
\frac{d\mu_\ba}{d\lam_1}(s_0) \ + \
\frac{d\mu_\ba}{d\lam_1}(b_0) \ + \
\frac{d\mu_\ba}{d\lam_1}(b_1) \ + \
\frac{d\mu_\ba}{d\lam_1}(b_1)
\quad\And\quad
  \dAX \ := \ \set{\ba\in\lAX}{M'(\ba)<\oo}.
\]
Clearly, $\dAX\subset\oAX$.}

The other two main results of this section are:
\Theorem{\label{threshold.convergence}}
{
  Fix $\gk\in\sK$.  Let $\{(s_0^n,b_0^n,b_1^n,s_1^n)\}_{n=1}^\oo\subset\Theta$,
with
$\D\lim_{n\goto\oo} (s_0^n,b_0^n,b_1^n,s_1^n)  =  (s_0,b_0,b_1,s_1)$.
For each $n\in\Natur$, let $\Life_n:\Ll\into\Ll$ be the RealLife EA
defined by $(s^n_0,b^n_0,b^n_1,s^n_1)$ and $\gk$. 
Then
\bthmlist
  \item  $\D\Lllim_{n\goto\oo} \Life_n(\ba) \ = \ \Life(\ba)$
for all $\ba\in\oAX$.

  \item  $\{\Life_n\}_{n=1}^\oo$ evolves to $\Life$ on $\oAX$.
\ethmlist
}

\Theorem{\label{akylea.convergence}}
{
  Fix $(s_0,b_0,b_1,s_1)\in\Theta$.
  Let $\{\gk_n\}_{n=1}^\oo \subset \sK$ be such that
$\D\Lllim_{n\goto\oo} \gk_n = \gk$.
For each $n\in\Natur$, let $\Life_n:\Ll\into\Ll$ be the RealLife EA
defined by $(s_0,b_0,b_1,s_1)$ and $\gk_n$. 
Then
\bthmlist
  \item  $\D\Lllim_{n\goto\oo} \Life_n(\ba) \ = \ \Life(\ba)$
for all $\ba\in\oAX$.

  \item  If  $\D\sup_{n\in\Natur} \ \norm{\gk_n}{\oo} \ < \ \oo$,  then
$\{\Life_n\}_{n=1}^\oo$ evolves to $\Life$ on $\oAX$.
\ethmlist
}

\Corollary{}
{  Fix $(s_0,b_0,b_1,s_1)\in\Theta$.
Let $\dK\subset\dX$ and $\gk:=\lam[\dK]^{-1}\chr{\dK}$.
\bthmlist
  \item
Let $\{\dK_n\subset\dX\}_{n=1}^\oo$ be such that
$\D \lim_{n\goto\oo} \lam[\dK\symdif\dK_n] = 0$.
For each $n\in\Natur$, let $\Life_n\colon\Ll\rightarrow\Ll$ be the RealLife EA
defined by $(s_0,b_0,b_1,s_1)$ and $\gk_n:=\lam[\dK_n]^{-1}\chr{\dK_n}$. 
Then  $\D\Lllim_{n\goto\oo} \Life_n(\ba)  =  \Life(\ba)$
for all $\ba\in\oAX$, and $\{\Life_n\}_{n=1}^\oo$ evolves
to $\Life$ on $\oAX$.

  \item Let  $\gG\in\sK$ be smooth. 
For any $n\in\Natur$, define $\gG_n\in\sK$ by
$\gG_n(x) := n^D \cdot \gG(n x)$, \  $\forall \ x\in\dX$, and then
let $\Life_n\colon\Ll\rightarrow\Ll$ be the RealLife EA
defined by $(s_0,b_0,b_1,s_1)$ and $\gk_n  :=  \gG_n*\gk$
 {\rm(a smooth kernel)}.
Then $\D\Lllim_{n\goto\oo} \Life_n(\ba)  =  \Life(\ba)$
for all $\ba\in\oAX$, and $\{\Life_n\}_{n=1}^\oo$ is evolves
to $\Life$ on $\oAX$.
\ethmlist
}
\bthmprf
{\bf(a)} \ 
Let $K:=\D\inf_{n\in\Natur} \, \lam[\dK_n]$, and
assume $K>0$, so $1/K < \oo$.  Now, $K\leq \lam[\dK]$, so
$\D \norm{\gk_n-\gk}{1}  <  (1/K)
  \lam[\dK_n\symdif\dK] \ \goesto{n\goto\oo}{} \ 0$, so 
Theorem \ref{akylea.convergence}(a) yields pointwise convergence
to $\Life$
on $\oAX$.
Also, $\D \sup_{n\in\Natur} \ \norm{\gk_n}{\oo} \ = \ 1/K <\oo$,
so Theorem \ref{akylea.convergence}(b) yield evolution.

{\bf(b)} \ 
The sequence $\{\gG_n\}_{n=1}^\oo\subset\sK$ is a convolutional
approximation of identity \cite[Thm.8.14(a)]{Folland}, so
$\D\Lllim_{n\goto\oo} \gk_n  =  \gk$, so
 Theorem \ref{akylea.convergence}(a)
yields pointwise convergence to
$\Life$ on $\oAX$.  Also, $\norm{\gk_n}{\oo} \leq
\norm{\gk}{\oo} <\oo$ by Young's inequality \cite[Thm.8.7]{Folland}, so 
Theorem \ref{akylea.convergence}(b) yields evolution.
\ethmprf

To prove Theorems \ref{reallife.is.continuous} to 
\ref{akylea.convergence},
we need some notation. For any $\del>0$, we define
\beq
\bW^s_\del & \ := \ &  \OO{s_0-\del,s_0+\del} \ \union \ 
\OO{s_1-\del,s_1+\del}\\
\And \ \bW^b_\del & \ := \ &  \OO{b_0-\del,b_0+\del} \ \union \
\OO{b_1-\del,b_1+\del}.
\eeq
If $\ba\in\AX$ and $\alp:=\gk*\ba$, then we define
$M^s_\ba,M^b_\ba:\OO{0,1}\into\Real$ by 
\[
\hspace{-2em}M^s_\ba(\del) \ :=\ \lam\lb[ \alp^{-1}(\bW^s_\del) \rb]
\quad \And\quad
M^b_\ba(\del) \ := \  \lam\lb[ \alp^{-1}(\bW^b_\del) \rb],
\quad\mbox{for all $\del> 0$.}
\]

\Lemma{\label{conv.lemma.2}}
{
Let $\ba,\ba'\in\Ll$.  Let $\alp:=\gk*\ba$ and $\alp':=\gk*\ba'$.
Then:
\bthmlist 
  \item $\norm{\Life(\ba)-\Life(\ba')}{1} 
\ \leq \ 2\cdot\norm{\ba-\ba'}{1} \ + \ 
\norm{\gs\circ\alp - \gs\circ\alp'}{1} \ + \
 \norm{\gb\circ\alp - \gb\circ\alp'}{1}$.

  \item $\norm{\gs\circ\alp' - \gs\circ\alp}{1} \ \leq \ 
M^s_\ba\lb(\norm{\alp'-\alp}{\oo}\rb)$
and $\norm{\gb\circ\alp' - \gb\circ\alp}{1} \ \leq \ 
M^b_\ba\lb(\norm{\alp'-\alp}{\oo}\rb)$.

  \item If $K:=\norm{\gk}{\oo}$, then
$\norm{\alp'-\alp}{\oo} \ \leq K\cdot \norm{\ba-\ba'}{1}$.

  \item  If $M(\ba)=0$, then 
$\D\lim_{\del\goto0} \ M^s_\ba(\del) \ = \ 0 
 \ = \ \lim_{\del\goto0} \ M^b_\ba(\del)$.
\ethmlist
}
\bthmprf {\bf(a)}  Eqn.(\ref{afterlife.CA2}) says
$\Life(\ba)  =  \ba\cdot \gs\circ\alp \ + \ (1-\ba)\cdot \gb\circ\alp$
and 
$\Life(\ba')  =  \ba'\cdot \gs\circ\alp' \ + \ (1-\ba')\cdot \gb\circ\alp'$.
Thus,
\beq
\hspace{-3em}
\Life(\ba)-\Life(\ba') &=&
(\ba-\ba')\cdot \gs\circ\alp  \ + \
\ba'\cdot \lb(\gs\circ\alp - \gs\circ\alp'\rb) \\
&&\qquad \ + \
(\ba'-\ba)\cdot \gb\circ\alp  \ + \
 (1-\ba')\cdot \lb(\gb\circ\alp - \gb\circ\alp'\rb). 
\\
\hspace{-3em}
\therefore \ 
\norm{\Life(\ba)-\Life(\ba')}{1} &\leq&
\norm{\ba-\ba'}{1}\cdot \norm{\gs\circ\alp}{\oo}  \ + \
\norm{\ba'}{\oo}\cdot \norm{\gs\circ\alp - \gs\circ\alp'}{1}
\\ \nonumber
 &&\quad
\ + \ \norm{\ba'-\ba}{1}\cdot \norm{\gb\circ\alp}{\oo} 
 \ + \
 \norm{1-\ba'}{\oo}\cdot \norm{\gb\circ\alp - \gb\circ\alp'}{1}. \\
&\eeequals{(*)} & \nonumber
2\cdot\norm{\ba-\ba'}{1} \ + \ 
\norm{\gs\circ\alp - \gs\circ\alp'}{1} \ + \
 \norm{\gb\circ\alp - \gb\circ\alp'}{1},
\eeq
where $(*)$ is because $\norm{\gs}{\oo}=\norm{\gb}{\oo}=\norm{\ba'}{\oo}
=\norm{1-\ba'}{\oo}=1$.

 {\bf(b)}\quad
$\norm{\gs\circ\alp' - \gs\circ\alp}{1}  =  \lam[\Del]$,
where $\Del:=\set{x\in\dX}{\gs\circ\alp(x)\neq\gs\circ\alp'(x)}$.
Let $\del:=\norm{\alp'-\alp}{\oo}$; then
$ \Del\ \subseteq\ \alp^{-1}(\bW^s_\del)$, so
$ \lam[\Del]
 \leq  
\lam\lb[ \alp^{-1}(\bW^s_\del) \rb]  =  M_\ba^s(\del)$.  The proof for $M^b_\ba$ is analogous.

{\bf(c)} \quad
$\norm{\alp' - \alp}{\oo}
\ = \ 
\norm{\gk*(\ba'-\ba)}{\oo}
\ \leq \
\norm{\gk}{\oo}\cdot\norm{\ba'-\ba}{1}$,
by  Young's inequality.

{\bf(d)} \quad  $\D  \lim_{\del\goto0} \
  \lb(M^s_\ba(\del)  +  M^b_\ba(\del)\rb) \ \eeequals{(*)} \  
\lam\lb[\alp^{-1}\{s_0,s_1\}\rb]  + \lam\lb[\alp^{-1}\{b_0,b_1\}\rb] \ = \ 
M(\ba) \ = \ 0$.   To see $(*)$, note that $\{s_0,s_1\} = \Intsct_{\del>0} 
\bW_s^\del$, so $\alp^{-1}\{s_0,s_1\} = \Intsct_{\del>0} 
\alp^{-1}(\bW_s^\del)$, so
$\lam\lb[\alp^{-1}\{s_0,s_1\}\rb] = \D
\lim_{\del\goto 0} \lam\lb[\alp^{-1}(\bW_s^\del)\rb] = \lim_{\del\goto0}
M^s_\ba(\del)$, by `continuity from above' for measures
\cite[Thm.1.8(d)]{Folland}.
Likewise, $\lam\lb[\alp^{-1}\{b_0,b_1\}\rb] =\D \lim_{\del\goto0}
M^b_\ba(\del)$.
 \ethmprf

\bthmprf[Proof of Theorem {\rm\ref{reallife.is.continuous}:}]
Let $\ba\in\oAX$ and $\ba'\in\Ll$, with $\norm{\ba-\ba'}{1}<\del$.  
If $\alp := \gk*\ba$ and $\alp':=\gk*\ba'$, then
\beq
\norm{\Life(\ba)-\Life(\ba')}{1} 
& \leeeq{(*)}&
2\cdot\norm{\ba-\ba'}{1} \ + \ 
\norm{\gs\circ\alp - \gs\circ\alp'}{1} \ + \
 \norm{\gb\circ\alp - \gb\circ\alp'}{1}
\\ &\leeeq{(\dagger)} & 
2\del \ + \  M^s_\ba(K\del) \ + \ M^b_\ba(K\del)
\qquad\mbox{(where $K:=\norm{\gk}{\oo}$).}
\eeq
Here, $(*)$ is by Lemma \ref{conv.lemma.2}(a), and  
$(\dagger)$ is by Lemma \ref{conv.lemma.2}(b,c).  
But $M(\ba)=0$, so
$M^s_\ba(K\del) \ + \ M^b_\ba(K\del) \goesto{\del\goto0}{} 0$, by
Lemma \ref{conv.lemma.2}(d).  
\ethmprf

\bthmprf[Proof of Theorem {\rm \ref{oAX.is.comeager}:}]
{\em $\!\!\shift{}$-invariant:}  Fix $v\in\dX$.
Let $\ba\in\lAX$. Let $\ba_v:=\shift{v}(\ba)$.  If $\alp=\gk*\ba$ and
$\alp_v=\gk*\ba_v$, then $\alp_v = \shift{v}(\alp)$.  Thus,
$\alp_v^{-1}\{s_0,b_0,b_1,s_1\}  =  
\shift{v}\lb[\alp^{-1}\{s_0,b_0,b_1,s_1\}\rb]$, so
 $M(\ba_v) = \lam\lb[\alp_v^{-1}\{s_0,b_0,b_1,s_1\}\rb] =
 \lam\lb[\alp^{-1}\{s_0,b_0,b_1,s_1\}\rb] = M(\ba)$.
Thus, $\ba\in\oAX$ iff $\ba_v\in\oAX$.

\noindent{\em  Dense G$\del$:} \ 
  For any $r\in\CC{0,1}$, and  $m\in\Natur$, let $\sC_m(r)
:=\set{\alp\in\Loo}{\lam\lb[\alp^{-1}\{r\}\rb] \geq \frac{1}{m}}$.

\Claim{\label{oAX.is.comeager.C1}  $\sC_m(r)$ is closed in $\Loo$.}
\bclaimprf
  Let $\{\alp_n\}_{n=1}^\oo \subset \sC_m(r)$, and let
$\D\Loolim_{n\goto\oo} \alp_n =\alp$.  We claim that 
$\lam\lb[\alp^{-1}\{r\}\rb] \geq \frac{1}{m}$, hence
$\alp\in\sC_m(r)$.
  To see this, fix $\eps>0$, and let $\bU_\eps := \OO{r-\eps,r+\eps}$.
Find $n\in\Natur$ such that $\norm{\alp_n-\alp}{\oo}<\eps$.
Then for any $x\in\dX$,
\[
\hspace{-2.5em}\statement{$x\in\alp_n^{-1}\{r\}$}
\iff
 \statement{$\alp_n(x)=r$}
\implies
 \statement{$\alp(x)\in\OO{r-\eps,r+\eps}$}
\iff\statement{$x\in\alp^{-1}(\bU_\eps)$}.
\]
Hence, $\alp_n^{-1}\{x\}\subset\alp^{-1}(\bU_\eps)$, which means
that $\lam\lb[\alp^{-1}(\bU_\eps)\rb]\geq \lam\lb[\alp_n^{-1}\{x\}\rb]
\geq 1/m$.

  This holds for any $\eps>0$.  Thus, 
$\lam\lb[\alp^{-1}\{x\}\rb] = \D\lim_{\eps\searrow 0} \lam\lb[\alp^{-1}(\bU_\eps)\rb] \geq 1/m$.
\eclaimprf
 
   Claim \ref{oAX.is.comeager.C1}
implies that $\sO_m(r) := \compl{\sC_m(r)}$ is open in $\Loo$.   Define $\gK:\lAX\ni \ba\mapsto (\gk*\ba)\in\Loo$;  then
$\gK$ is continuous (by Lemma \ref{conv.lemma.2}(c)), so
$\gO_m(r):=\gK^{-1}\lb[\sO_m(r)\rb]$ is open in $\lAX$.

  It remains to show that  $\gO_m(r)$ is $\bL^1$-dense in $\lAX$.  To
do this, let $\ba\in\compl{\gO_m(r)}\subset\lAX$.
 hence, if $\alp:=\gk*\ba$, then $\lam[\alp^{-1}\{r\}]\geq 1/m$. 
 Let  $A:=\norm{\ba}{1}$ and $K:=\norm{\gk}{\oo}$.  We define constants
$L:=8(K+A)$ and $J:=(1+4LA/r)$.

\Claim{\label{oAX.is.comeager.C0}
For any $\eps>0$, there is some $\ba'\in\lAX$ such that
$\norm{\ba-\ba'}{1}<J\cdot\eps$ and so that, if $\alp':=\gk*\ba'$, then
$(\alp')^{-1}\{r\} \ \subseteq \ \alp^{-1}\OO{r+L\eps,r+3L\eps}$.}
\bclaimprf
Assume $\eps<A$. 
 Suppose $\ba=\chr{\bA}$ for some measurable $\bA\subset\dX$.  Let
$\bC := \bC_1\disj\bC_2\disj\cdots\disj\bC_N$ be a finite disjoint union
of open cubes such that $\lam[\bC\symdif\bA]<\eps$ 
\cite[Thm. 2.40(c), p.68]{Folland}.
Now, $\gk$ is almost continuous, so there is a thin set $\bT\subset\dX$ 
and  $\gam>0$ so that 
$\lam[\dB(\bT,2\gam)] < \eps$ and
so that $\gk$ is uniformly continuous on $\bY:=\compl{\dB(\bT,\gam)}$.

  For any $x,y\in\dX$ and $\del>0$, we'll write ``$x \closeto{\del} y$''
to mean $|x-y|<\del$.
Find $\del$ so that, for any $y,y'\in\bY$, \ $\statement{$y\closeto{\del} y'$}
\implies\statement{$\gk(y)\closeto{\eps}\gk(y')$}$.
Assume $\del<\gam$.   By subdividing 
the cubes $\{\bC_n\}_{n=1}^N$ if necessary, we can assume all cubes have
diameter less than $\del$. 

\ignore{ Note that
\[
A \quad =\quad \norm{\ba}{1}\quad=\quad
\lam[\bA] \quad \closeto{\eps} \quad 
\lam[\bC]\quad = \quad\sum_{n=1}^N \lam[\bC_n]
\quad \leq \quad N\cdot\del^D;
\]
hence $(A-\eps)/\del^D\leq N$.}
 
  Fix $x\in\dX$.  Let $\bT_x=x-\bT$.  By reordering $\{\bC_n\}_{n=1}^N$
if necessary, we can find some $M_x < N$ such that,
$\bC_m\intsct \dB(\bT_x,\gam)  \neq  \emptyset$ for all $m\in\CC{1...M_x}$,
and
$\bC_n\intsct \dB(\bT_x,\gam)   =  \emptyset$ for all $n\in\OC{M_x...N}$.
Thus, if we define  $\bC_x  :=   \Disj_{m=1}^{M_x} \bC_m$, then
$\dB(\bT_x,\gam)  \subseteq \bC_x$.

\subclaim{\label{oAX.is.comeager.C0.1}
$\D\lam[\bC_x]<\eps$ for all $x\in\dX$.}
\bsubclaimprf
If $m\in\CC{1...M_x}$, then $
\bC_m \ \suuubset{(*)}\  \dB(\bT_x,\gam+\del) \ \suuubset{(\dagger)}\
  \dB(\bT_x,2\gam)$.
Here, $(*)$ is because $\dB(\bT_x,\gam)\intsct\bC_m\neq0$ and $\diam{\bC_m}<\del$.
$(\dagger)$ is because $\del<\gam$.
  Thus, 
$\bC_x \subseteq \dB(\bT_x,2\gam)$, so $\lam[\bC_x]\leq \lam[\dB(\bT_x,2\gam)]
< \eps$ by definition of $\gam$.  This works for all $x\in\dX$.
\esubclaimprf
 For each $n\in\OC{M_x...N}$, fix some $c_n\in\bC_n$,
and let $k_n := \gk(x-c_n)$.  

\subclaim{\label{oAX.is.comeager.e2a} 
For any $n\in\OC{M_x...N}$ and every $c\in\bC_n$, \ 
$\gk(x-c)\closeto{\eps}k_n$.}
\bclaimprf
  $c\not\in\dB(\bT_x,\gam)$, so
 $(x-c)\not\in\dB(\bT,\gam)$, so  $(x-c)\in\bY$.
  Likewise, $(x-c_n)\in\bY$. If $c\in\bC_n$, then
 $c\closeto{\del}c_n$, so $(x-c)\closeto{\del}(x-c_n)$, so
$\gk(x-c)\closeto{\eps}k_n$, by definition of $\del$.
\eclaimprf

\subclaim{\label{oAX.is.comeager.e2} 
$\D \gk*\ba(x) \ \ \closeto{L\eps/4} \ \sum_{n=M_x+1}^N k_n \lam[\bC_n]$.}
\bsubclaimprf
Claim \ref{oAX.is.comeager.C0}.\ref{oAX.is.comeager.e2a} implies that,
for any $n\in\OC{M_x...N}$,
\[
 \lb| \int_{\bC_n}  \gk(x-c) \ dc \ - \ k_n\cdot\lam[\bC_n] \rb|
\quad \leq \quad \eps \cdot\lam[\bC_n].\]
\begin{eqnarray}
\hspace{-2em}\mbox{Thus,}\quad
\nonumber
\lefteqn{\lb|\sum_{n=M_x+1}^N \int_{\bC_n}  \gk(x-c) \ dc \ - \ 
\sum_{n=M_x+1}^N  k_n\lam[\bC_n] \rb|}
\\\hspace{-3em} & \leq &
 \eps \cdot\sum_{n=M_x+1}^N  \lam[\bC_n] 
\quad\leq\quad \eps\cdot\lam[\bC] 
\label{oAX.is.comeager.e1}
\quad \leeeq{(*)} \quad
 \eps\cdot(A+\eps) \quad\leeeq{(\dagger)}\quad 2A\eps.
\qquad\qquad
\end{eqnarray}
Here, $(*)$ is because $\lam[\bA\symdif\bC]<\eps$,
and $(\dagger)$ is because $\eps<A$. Thus, if $K:=\norm{\gk}{\oo}$, then
\beq
\hspace{-2em}
\lefteqn{\gk*\ba(x) \quad =\quad  \int_\dX \gk(x-y)\cdot\ba(y) \ dy
\quad=\quad \int_\bA \gk(x-y) \ dy
\quad\stackrel{^{(*)}}{\closeto{K\eps}}\quad
 \int_\bC \gk(x-c) \ dc}
\\ \hspace{-2em}&=& 
\sum_{n=1}^N \int_{\bC_n} \gk(x-c) \ dc
\quad=\quad
\sum_{m=1}^{M_x} \int_{\bC_m} \gk(x-c) \ dc \ + \
\sum_{n=M_x+1}^N \int_{\bC_n} \gk(x-c) \ dc \\
\hspace{-2em}&\stackrel{^{(\dagger)}}{\closeto{2A\eps}}&
 \int_{\bC_x} \gk(x-c) \ dc   
\ + \ \sum_{n=M_x+1}^N k_n \lam[\bC_n] 
\quad\stackrel{^{(\ddagger)}}{\closeto{K \eps}} \quad \sum_{n=M_x+1}^N k_n \lam[\bC_n]. 
\eeq
Here, $(*)$ is because $\lam[\bA\symdif\bC]<\eps$, 
\ $(\dagger)$ is by \eref{oAX.is.comeager.e1}, and
$(\ddagger)$ is by Claim \ref{oAX.is.comeager.C0}.\ref{oAX.is.comeager.C0.1}.
The claim follows because $K\eps + 2A\eps + K\eps \ = \ 2(K+A)\eps \ = \ 
L\eps/4$ [because $L:=8(K+A)$].
\esubclaimprf
If $\eps$ is small enough, then $2L\eps/r<1$.  Thus
if $\beta:=(1+2L\eps/r)^{-1}$, then $\frac{1}{2}<\beta<1$.
For each $n\in\CC{1...N}$, let $\bC'_n$ be
the cube obtained by multiplying $\bC_n$ by $\sqrt[D]{\beta}$ but keeping 
the centre the same.  Hence $\lam[\bC'_n] =\beta \cdot \lam[\bC_n]$.
Let $\bA':=\Disj_{n=1}^N \bC'_n$  and let $\ba':=\chr{\bC'}$.  Then
\beq
  \norm{\ba-\ba'}{1} &=&\quad
\lam[\bA\symdif\bA']
\ \ \leq\ \ 
\lam[\bA\symdif\bC]
+
\lam[\bC\symdif\bA']
\\&\leq&\quad
\eps \ + \ (1-\beta)\cdot\lam[\bC]
\ \ \leeeq{(*)}\ \ 
(1+4LA/r)\eps
\ = \ J\eps.
\eeq
Here, $(*)$ is because  $\lam[\bC]<2\lam[\bA]=2A$, and 
$\D 1-\beta \ = \ 1 - \frac{1}{1+2L\eps/r} \ = \ 
\frac{2L\eps/r}{1+2L\eps/r} \ < \  2L\eps/r$.

\subclaim{\label{oAX.is.comeager.C0.2}
Let $\alp:=\gk*\ba$ and $\alp':=\gk*\ba'$.
Then
$\norm{\alp' - \beta\cdot\alp}{\oo} < L\eps/2$.}
\bsubclaimprf 
For each $x\in\dX$, Claim \ref{oAX.is.comeager.C0}.\ref{oAX.is.comeager.e2}
yields some
suitable reordering of $\{\bC_n\}_{n=1}^N$ and some suitable 
$M_x\in\CC{1...N}$ such that
\beq
\hspace{-1em}
\alp'(x)&=&
\gk*\ba'(x)
\ \ \stackrel{^{(*)}}{\closeto{\frac{L\eps}{4}}} \sum_{n=M_x+1}^N k_n \lam[\bC'_n] 
\ \ =\ \  \beta \sum_{n=M_x+1}^N k_n \lam[\bC_n]  
\\
& \stackrel{^{(\dagger)}}{\closeto{\beta\frac{L\eps}{4}}} & \quad \beta \bk*\ba(x) 
\ \ =\ \ \beta\cdot\alp(x).
\eeq
  Here,
$(\dagger)$ is by 
Claim \ref{oAX.is.comeager.C0}.\ref{oAX.is.comeager.e2}, and $(*)$ is by
an argument identical to Claim
 \ref{oAX.is.comeager.C0}.\ref{oAX.is.comeager.e2}. 
Thus, for all $x\in\dX$, we have
 $|\alp'(x)-\beta\cdot\alp(x)| < \frac{L\eps}{4} + \beta\frac{L\eps}{4} < L\eps/2$.
 The claim follows.
\esubclaimprf
  
  Let $r':= r/\beta \ = \ r(1+2L\eps/r) \ = \ r+2L\eps$. 
Then 
\[
(\alp')^{-1}\{r\} \quad \suuubset{(*)} \quad \alp^{-1}\OO{r'-L\eps,r'+L\eps}
\quad \eeequals{(\dagger)} \quad \alp^{-1}\OO{r+L\eps,r+3L\eps}.
\]
Here, $(\dagger)$ is
 because $r'-L\eps =  r + 2L\eps - L\eps
= r-L\eps$, and $r'+L\eps =  r+3L\eps$.  $(*)$ is because \
$\statement{$x\in(\alp')^{-1}\{r\}$}
\implies
\statement{$\bet \alp(x) \stackrel{(\sharp)}{\closeto{\frac{L\eps}{2}}}\alp'(x) = r$}
\implies 
\statement{$\alp(x) \closeto{\frac{L\eps}{2\bet}} r/\bet = r'$}
\iiimplies{(\flat)} 
\statement{$\alp(x) \closeto{L\eps} r'$}$,
\ 
where $(\sharp)$ is by
Claim \ref{oAX.is.comeager.C0}.\ref{oAX.is.comeager.C0.2},
and $(\flat)$ is because  $L\eps/2\bet < L\eps$
because $\bet>1/2$.
\eclaimprf

\Claim{\label{oAX.is.comeager.C2} 
For any $r\in\CC{0,1}$, and any $m\in\Natur$, \ 
$\gO_m(r)$ is $\bL^1$-dense in $\lAX$.}
\bclaimprf
  Let $\ba\in\compl{\gO_m(r)}\subset\lAX$.
  Fix $\eps>0$;  we want some
$\ba'\in\gO_m(r)$ with $\norm{\ba-\ba'}{1}<\eps$. 

 Let $L$ and $J$ be the constants defined prior to Claim \ref{oAX.is.comeager.C0}.
 For each $n\in\Natur$, let $\eps_n:=\eps/(L 4^n)$, and  apply Claim 
\ref{oAX.is.comeager.C0} to obtain a sequence $\{\ba_n\}_{n=1}^\oo$ such that
$\norm{\ba_n-\ba}{1} < J \eps_n = J\eps/(L 4^n)$, and so that,
if $\alp_n:=\gk*\ba_n$, then 
$\alp_n^{-1}\{r\} \subset \alp^{-1}\OO{r+ \frac{1}{4^n}\eps,
r+\frac{3}{4^n}\eps}$.

\subclaim{\label{oAX.is.comeager.C2.1}
There are infinitely many $n\in\Natur$ such that
 $\ba_n\in\gO_m(r)$.}
\bsubclaimprf
Suppose not.  Then there is some $N\in\Natur$ so that
 $\lam[\alp_n^{-1}\{r\}]\geq\frac{1}{m}$ for all $n\geq N$.  This means that 
$\lam\lb[\alp^{-1}\OO{r+ \frac{\eps}{4^n},
r+\frac{3\eps}{4^n}}\rb] \geq\frac{1}{m}$.
But the sequence of open intervals
\[
\hspace{-2em}
\OO{r+ \mbox{$\frac{\eps}{4}$},
r+\mbox{$\frac{3\eps}{4}$}}, 
\ 
\OO{r+ \mbox{$\frac{\eps}{16}$},
r+\mbox{$\frac{3\eps}{16}$}}, 
\ 
\OO{r+ \mbox{$\frac{\eps}{64}$},
r+\mbox{$\frac{3\eps}{64}$}}, \
\ldots, \ 
\OO{r+ \mbox{$\frac{\eps}{4^n}$},r+\mbox{$\frac{3\eps}{4^n}$}},\ldots
\] 
are disjoint.  Hence, the sequence of subsets
\[
\hspace{-1em}
\alp^{-1}\OO{r+ \mbox{$\frac{\eps}{4}$},
r+\mbox{$\frac{3\eps}{4}$}}, 
\quad 
\alp^{-1}\OO{r+ \mbox{$\frac{\eps}{16}$},
r+\mbox{$\frac{3\eps}{16}$}}, 
\quad 
\alp^{-1}\OO{r+ \mbox{$\frac{\eps}{64}$},
r+\mbox{$\frac{3\eps}{64}$}}, \
\ldots, 
\] 
of also disjoint.  All of these are subsets of $\alp^{-1}\OO{r,r+\eps}$,
which means that 
\[
\hspace{-1em}
\lam\lb[\alp^{-1}\OO{r,r+\eps}\rb] 
\ \ \geq \ \ \sum_{n=N}^\oo \lam\lb[\alp^{-1}\OO{r+ \mbox{$\frac{\eps}{4^n}$},
r+\mbox{$\frac{3\eps}{4^n}$}}\rb]
\ \ \geq \ \  \sum_{n=N}^\oo \frac{1}{m} \ \ = \ \  \oo.
\]
  But $\lam\lb[\alp^{-1}\OO{r,r+\eps}\rb]   \stackrel{\scriptscriptstyle(*)}{\leq}  \frac{1}{r}
\norm{\alp}{1}  \stackrel{\scriptscriptstyle(\dagger)}{\leq}  \frac{1}{r}
\norm{\gk}{\oo}\cdot\norm{\ba}{1}  =  K\cdot A/r$, which is finite.
Contradiction.  Here, $(*)$ is Chebyshev's inequality
\cite[Thm.6.17]{Folland}, and $(\dagger)$ is Young's
inequality.  \esubclaimprf

Let $\ba_n\in\gO_m(r)$ be as in 
Claim \ref{oAX.is.comeager.C2}.\ref{oAX.is.comeager.C2.1},
with $n$ large enough that $J/L 4^n<1$. 
Thus,
$\ba_n\closeto{\eps}\ba$. Since $\eps$ was arbitrary, we conclude
that $\ba$ is a cluster point of $\gO_m(r)$.  Since $\ba\in\compl{\gO_m(r)}$
was arbitrary, we conclude that $\gO_m(r)$ must be dense in $\lAX$.
\eclaimprf
For any $m\in\Natur$,
 Claim \ref{oAX.is.comeager.C2} means that
 $\gO_m:=\gO_m(s_0)\intsct\gO_m(b_0)\intsct\gO_m(b_1)\intsct\gO_m(s_1)$
is open and dense in $\lAX$.  Thus, 
$\oAX = \Intsct_{m=1}^\oo \gO_m$ is dense $G\delta$.
\ethmprf

{\em Remark:}  Note that the proofs of shift-invariance and `$G\delta$'
in Theorem \ref{oAX.is.comeager} 
do not depend on the almost-continuity of $\gk$.\hfill$\diamondsuit$

\breath

Now we'll prove Theorems \ref{threshold.convergence} and
\ref{akylea.convergence}.  Let
$\{\Phi_n:\Ll\into\Ll\}_{n=1}^\oo$ be a family of transformations of
$\Ll$.  If $\ba\in\Ll$, then $\{\Phi_n\}_{n=1}^\oo$ is 
{\dfn eventually equicontinuous (EE)} at
$\ba$ if, for any $\gam>0$, there is some $\del>0$ and some
$N\in\Natur$ so that, for any $\ba'\in\Ll$, 
if $\ba'\closeto{\del}\ba$, then for all $n\geq N$, \
$\Phi_n(\ba')\closeto{\gam}\Phi_n(\ba)$.  The sequence
$\{\Phi_n\}_{n=1}^\oo$ is {\dfn equicontinuous} if it is EE, and we
can further hold $N:=1$ for all $\gam>0$. 

\Proposition{\label{fixed.point.convergence.lemma}}
{
 Let $\Phi:\Ll\into\Ll$
and $\{\Phi_n\colon\Ll\into\Ll\}_{n=1}^\oo$ be a sequence of operators,
and let $\gA\subset\Ll$ be a $\shift{}$-invariant subset.
Suppose $\{\Phi_n\}_{n=1}^\oo$  is EE
at all points in $\gA$, and that $\{\Phi_n\}_{n=1}^\oo$
converges to $\Phi$ pointwise on $\gA$. Then
$\{\Phi_n\}_{n=1}^\oo$ evolves to $\Phi$ on $\gA$.
}
\bthmprf 
 Let $\ba\in\gA$; then
$\D\lim_{n\goto\oo} \Phi_n(\ba)  =  \Phi(\ba)$.
Let $\{\ba_n\}_{n=1}^\oo\subset\Ll$ be such that
 $\D\lim_{n\goto\oo} \ba_n  =  \ba$.
We will verify parts {\bf(a)} and {\bf(b)} in the definition of `evolution'.

{\bf(a)}\quad Suppose $\Phi_n(\ba_n)=\ba_n$ for all $n\in\Natur$. 
We claim that $\Phi(\ba) \ = \ \ba$.

To see this, fix $\gam>0$.  By eventual equicontinuity, 
find $\del>0$ and $N\in\Natur$
so that,
if $\norm{\ba -  \ \ba'}{1}<\del$, then
$\norm{\Phi_n(\ba) -  \ \Phi_n(\ba')}{1}<\frac{\gam}{3}$ for all $n\geq N$.
Assume that $\del<\frac{\gam}{3}$.  If $n\geq N$ is large enough, then
$\norm{\ba -  \ \ba_n}{1}<\del<\frac{\gam}{3}$, and also
$\norm{\Phi(\ba) -  \ \Phi_n(\ba)}{1}<\frac{\gam}{3}$.   Thus,
\beq
\hspace{-2em} \norm{\Phi(\ba) -  \ \ba}{1}
&\leq & \ \ 
\norm{\Phi(\ba) -  \ \Phi_n(\ba)}{1}
\ + \ 
\norm{\Phi_n(\ba) -  \ \Phi_n(\ba_n)}{1}
\ + \ 
\norm{\Phi_n(\ba_n) -  \ \ba}{1} \\
\hspace{-2em}
&\leq& \ \ 
\frac{\gam}{3}
\ + \ 
\frac{\gam}{3}
\ + \ 
\norm{\ba_n -  \ \ba}{1} 
\quad\leq\quad
\frac{2\gam}{3} + \del \quad\leq\quad\gam.
\eeq
Since $\gam$ is arbitrary, we conclude that $\Phi(\ba)=\ba$.

{\bf(b)} {\bf[i]} is simply a special case of {\bf[ii]} with $v:=0$, so we will
prove {\bf[ii]}.  Let $\ba\in\gA$, and
suppose, for each $p\in\CC{0...P}$, that $\Phi^p(\ba)\in\gA$.
Define $\Psi:=\shift{-v}\circ\Phi$.
 For all $p\in\CC{0...P}$, let 
$\ba^p:=\Psi^p(\ba)=\shift{-pv}\circ\Phi^p(\ba)$.  
Then  $\ba^p\in\gA$, because $\gA$ is $\shift{}$-invariant.

For all $n\in\Natur$, let $\Psi_n := \shift{-v_n}\circ\Phi_n$.
Thus, for any $p\in\CC{0...P}$, \  $\Psi^p_n = \shift{-pv_n}\circ\Phi^p_n$.
Suppose that $\Phi^P_n(\ba_n) = \shift{Pv_n}(\ba_n)$
for all $n\in\Natur$ [or equivalently, that $\Psi^P_n(\ba_n)=\ba_n$].
We claim that $\Phi^P(\ba) = \shift{Pv}(\ba)$
[or equivalently, that $\Psi^P(\ba)=\ba$].

\Claim{\label{fixed.point.convergence.lemma.C1}
 For all $q\in\CC{0...P}$, the sequence $\{\Psi_n\}_{n=1}^\oo$
is eventually equicontinuous at $\ba^q$.}
\bclaimprf
 $\{\Phi_n\}_{n=1}^\oo$ is EE at $\ba^q$
because $\ba^q\in\gA$.
Thus $\{\Psi_n\}_{n=1}^\oo$ is
EE at $\ba^q$, because $\Psi_n = \shift{-v_n}\circ\Phi_n$, and
the $\bL^1$ metric is $\shift{}$-invariant.
\eclaimprf

\Claim{\label{fixed.point.convergence.lemma.C2}
For all $p\in\CC{1...P}$, \ 
$\D\lim_{n\goto\oo} \Psi_n^p(\ba) \ = \ \Psi^p(\ba) \ = \ \ba^{p}$.
}
\bclaimprf[Proof {\rm(}by induction on $p$\rm{):}]
 For all $q\in\CO{0...P}$, \  $\ba^q\in\gA$, so $\D
 \lim_{n\goto\oo} \Phi_n(\ba^q)  =  \Phi(\ba^q)$.  Thus,
\begin{eqnarray}
\nonumber
 \lim_{n\goto\oo} \Psi_n(\ba^q) &=&
\lim_{n\goto\oo} \shift{-v_n}\circ\Phi_n(\ba^q) 
\ \   \eeequals{(*)} \ \ 
\lim_{n\goto\oo} \shift{-v_n}\circ\Phi(\ba^q) 
\\&=&
 \shift{-v}\circ\Phi(\ba^q) \ \  = \ \  \ba^{q+1}.
\label{life.fixed.point.convergence.e1}
\end{eqnarray}
Here, $(*)$ is because $\lim_{n\goto\oo} \Phi_n(\ba^q)  =  \Phi(\ba^q)$,
and the $\bL^1$ metric is $\shift{}$-invariant.
To obtain the base case $p=1$, set $q:=0$ in
eqn.(\ref{life.fixed.point.convergence.e1}). 

Now, let $q\in\CO{1...P}$, and suppose the claim is true for $q$.  Let $p:=q+1$.
Then
\[
\D \lim_{n\goto\oo} \Psi_n^p(\ba)
\quad=\quad
\lim_{n\goto\oo} \Psi_n \circ \Psi_n^q(\ba)
\quad\eeequals{(*)}\quad \lim_{n\goto\oo} \Psi_n (\ba^q)
\quad\eeequals{(\dagger)}\quad  \ba^{q+1} \ = \ \ba^p.
\]
Here, $(*)$ is because $\D\lim_{n\goto\oo} \Psi_n^q(\ba) 
 =  \ba^{q}$ by induction,
and because $\{\Psi_n\}_{n=1}^\oo$ is
EE at $\ba^q$ by Claim~\ref{fixed.point.convergence.lemma.C1}.
\  $(\dagger)$ is by eqn.(\ref{life.fixed.point.convergence.e1}).
\eclaimprf

\Claim{\label{fixed.point.convergence.lemma.C3}
For all $p\in\CC{1...P}$, \ 
 The sequence $\{\Psi^p_n\}_{n=1}^\oo$ is eventually equicontinuous at
 $\ba$.
}
\bclaimprf[Proof {\rm(}by induction on $p$\rm{):}]
  The base case is  Claim 
\ref{fixed.point.convergence.lemma.C1}.
Let $q\in\CO{0...P}$, and suppose  $\{\Psi^q_n\}_{n=1}^\oo$ is EE
at $\ba$.  Let $p:=q+1$.
 Fix $\gam>0$.  Claim 
\ref{fixed.point.convergence.lemma.C1}
yields $\del_1>0$ and $N_0\in\Natur$ such that:
\beqn
\label{life.fixed.point.convergence.e2}
\mbox{For all $n\geq N_0$,}\quad
\statement{$\ba'\closeto{2\del_1}\ba^q$}
\implies
\statement{$\Psi_n(\ba')\closeto{\gam/3}\Psi_n(\ba^q)$}.
\eeqn
Claim 
\ref{fixed.point.convergence.lemma.C2} 
says $\D\lim_{n\goto\oo} \Psi^{p}_n(\ba) = \ba^{p}$
and 
$\D\lim_{n\goto\oo} \Psi^q_n(\ba) = \ba^{q}$.
The same argument can be adapted to show
$\D\lim_{n\goto\oo} \Psi_n(\ba^q) = \Psi(\ba^q) = \ba^{q+1} = \ba^p$.
Thus, find $N_1\geq N_0$ so that, for all $n\geq N_1$,
\beqn
\label{life.fixed.point.convergence.e3}
\hspace{-2.5em}
\ {\bf(a)} \ 
\Psi^{p}_n(\ba) \  \closeto{\gam/3}  \ \ba^{p},
\qquad
{\bf(b)} \ 
\Psi^q_n(\ba) \  \closeto{\del_1} \ \ba^{q} 
\ \
\And \ \
{\bf(c)} \ 
\Psi_n(\ba^q)\  \closeto{\gam/3} \ \ba^{p} 
\qquad\eeqn
  By induction,  $\{\Psi^q_n\}_{n=1}^\oo$ is EE
at $\ba$, so find $\del>0$ and $N\geq N_1$ such
that,
\beq
\hspace{-2em}
\lefteqn{\mbox{For all $n\geq N$,}
\quad
\statement{$\ba'\closeto{\del}\ba$}
\quad \iiimplies{(EE)}\quad
\statement{$\Psi^q_n(\ba')\closeto{\del_1}\Psi^q_n(\ba)$}
\quad\iiimplies{(\diamond)}\quad
\statement{$\Psi^q_n(\ba') \closeto{2\del_1} \ba^q$}} \\
\hspace{-2em}
&\implies &
\statement{$\Psi^{p}_n(\ba') \ =  \ \Psi_n(\Psi^{q}_n(\ba')) 
\quad\stackrel{(*)}{\closeto{\gam/3}}\quad \Psi_n(\ba^q)
\quad\stackrel{(\dagger)}{\closeto{\gam/3}}\quad\ba^{p}
\quad\stackrel{(\ddagger)}{\closeto{\gam/3}}\quad \Psi^{p}_n(\ba)$}
\\ \hspace{-2em} &\implies&
\statement{$\Psi^{p}_n(\ba') \ \closeto{\gam} \ \Psi^{p}_n(\ba)$}.
\eeq  
Here, $(\diamond)$ is  by eqn.(\ref{life.fixed.point.convergence.e3}b),
\ $(*)$ is by \eref{life.fixed.point.convergence.e2},
\ $(\dagger)$ is by eqn.(\ref{life.fixed.point.convergence.e3}c),
and $(\ddagger)$ is  by eqn.(\ref{life.fixed.point.convergence.e3}a).
Thus,  $\{\Psi^p_n\}_{n=1}^\oo$ is EE.
\eclaimprf
Set $p:=P$ in Claim \ref{fixed.point.convergence.lemma.C3}  to conclude
$\{\Psi^P_n\}_{n=1}^\oo$ is EE at $\ba$.
Set $p:=P$ in Claim \ref{fixed.point.convergence.lemma.C2}  to get
$\D\lim_{n\goto\oo} \Psi^P_n(\ba)  =  \Psi^P(\ba)$.
By hypothesis, $\Psi^P_n(\ba_n)=\ba_n$ for all $n\in\Natur$.
Now apply part {\bf(a)} to the sequence $\{\Psi^P_n\}_{n=1}^\oo$ to
conclude  $\Psi^P(\ba)  =  \ba$; hence $\Phi^P(\ba)=\shift{Pv}(\ba)$,
as desired.
\ethmprf
{\em Remarks:} Proposition \ref{fixed.point.convergence.lemma} doesn't
 need the functions $\{\Phi_n\}_{n=1}^\oo$ to be continuous
anywhere except at $\ba$, nor to converge to $\Phi$ anywhere except at
$\ba$.  Also, $\Ll$ could be replaced with any space $\gL$
of functions on $\dX$ equipped with a $\shift{}$-invariant
metric $d$ such that $\D\lim_{v\goto0} d(\shift{v}(\ba),\ba)=0$
for any $\ba\in\gL$.

\bthmprf[Proof of Theorem {\rm\ref{akylea.convergence}:}]
{\bf(a)} Fix $\ba\in\oAX$.  Let $\alp=\gk*\ba$,
and for all $n\in\Natur$, let $\alp_n:=\gk_n*\ba$.
Eqn.(\ref{afterlife.CA2}) says
$\Life(\ba) \ = \ \ba\cdot \gs\circ\alp \ + \ (1-\ba)\cdot \gb\circ\alp$
and
$\Life_n(\ba) \ = \ \ba\cdot \gs\circ\alp_n \ + \ (1-\ba)\cdot \gb\circ\alp_n$.
\beq
\hspace{-3em} \mbox{Thus,} \quad
\Life_n(\ba)-\Life(\ba) &=&
\ba\cdot \lb(\gs\circ\alp_n - \gs\circ\alp\rb) 
\ + \
 (1-\ba)\cdot \lb(\gb\circ\alp_n - \gb\circ\alp\rb), \\
\hspace{-3em} \mbox{so}\quad
\norm{\Life_n(\ba)-\Life(\ba)}{1} &\leq&
\norm{\ba}{\oo}\cdot \norm{\gs\circ\alp_n - \gs\circ\alp}{1}
\ + \
 \norm{1-\ba}{\oo}\cdot \norm{\gb\circ\alp_n - \gb\circ\alp}{1} \\
&\leq &
\norm{\gs\circ\alp_n - \gs\circ\alp}{1} \ + \
 \norm{\gb\circ\alp_n - \gb\circ\alp}{1}\\
&\leeeq{(*)} & 
M^s_\ba\lb(\norm{\alp_n-\alp}{\oo}\rb) \ + \ 
M^b_\ba\lb(\norm{\alp_n-\alp}{\oo}\rb),
\eeq
where $(*)$ is by  Lemma \ref{conv.lemma.2}(b).  But
\[
\hspace{-2em}
\norm{\alp_n-\alp}{\oo}
\  =  \ \norm{(\gk_n-\gk)*\ba}{\oo}
 \ \ \leeeq{(*)} \ \ 
 \norm{\gk_n-\gk}{1}\cdot\norm{\ba}{\oo}
\   =  \ \norm{\gk_n-\gk}{1} \  \goesto{n\goto\oo}{}  \ 0.
\]
$(*)$ is by Young's inequality.
Thus, $\norm{\Life_n(\ba)-\Life(\ba)}{1} \ \goesto{n\goto\oo}{}  \ 0$,
by Lemma \ref{conv.lemma.2}(d).

{\bf(b)}  The set $\oAX$ is $\shift{}$-invariant.
Given part {\bf(a)} together with
Proposition \ref{fixed.point.convergence.lemma}, it suffices to
show that  $\{\Life_n\}_{n=1}^\oo$ is $\Ll$-equicontinuous at
every $\ba\in\oAX$.

To see this, fix $\eps>0$.  Let $K:=\D\sup_{n\in\Natur} \ \norm{\gk_n}{\oo}$.
Let $\ba'\in\Ll$, with $\norm{\ba-\ba'}{1}<\del$.
Then Lemma  \ref{conv.lemma.2}(a,b,c) implies that:
\beqn
\label{akylea.convergence.e1}
\hspace{-2em}
\mbox{$\forall \ n\in\Natur$,} \qquad
\norm{\Life_n(\ba)-\Life_n(\ba')}{1} \quad\leq\quad
2\del \ + \  M^s_\ba(K\del) \ + \ M^b_\ba(K\del).
\qquad
\eeqn
Now,  $M(\ba)=0$, so  Lemma \ref{conv.lemma.2}(d) 
yields some $\del$ such
that  $M^s_\ba(K\del)  +  M^b_\ba(K\del)<\eps/2$.
Assume  $\del<\eps/4$. Then \eref{akylea.convergence.e1}
says that $\norm{\Life_n(\ba)-\Life_n(\ba')}{1} \ < \ \eps$ for
all $n\in\Natur$.
\ethmprf

\bthmprf[Proof of Theorem {\rm\ref{threshold.convergence}:}]
{\bf(a)} \ 
Fix $\ba\in\oAX$.  Let $\alp=\gk*\ba$.
For all $n\in\Natur$, let $\gs_n:=\chr{\CC{s^n_0,s^n_1}}$
and $\gb_n:=\chr{\CC{b^n_0,b^n_1}}$.  
Then
$\Life(\ba) \ = \ \ba\cdot \gs\circ\alp \ + \ (1-\ba)\cdot \gb\circ\alp$
and
$\Life_n(\ba) \ = \ \ba\cdot \gs_n\circ\alp \ + \ (1-\ba)\cdot \gb_n\circ\alp$.
Thus,
\begin{eqnarray}
\nonumber
\lefteqn{\hspace{-2em} \norm{\Life_n(\ba)-\Life(\ba)}{1}
\quad=\quad
\norm{\ba\cdot \lb(\gs_n\circ\alp - \gs\circ\alp\rb) 
\ + \
 (1-\ba)\cdot \lb(\gb_n\circ\alp - \gb\circ\alp\rb)}{1}} \\
\nonumber
 &\leq& \
\norm{\ba}{\oo}\cdot \norm{\gs_n\circ\alp - \gs\circ\alp}{1}
\ + \
 \norm{1-\ba}{\oo}\cdot \norm{\gb_n\circ\alp - \gb\circ\alp}{1} \\
&\leq & \
\norm{\gs_n\circ\alp - \gs\circ\alp}{1} \ + \
 \norm{\gb_n\circ\alp - \gb\circ\alp}{1}
\nonumber \\&=& \
\norm{\gs'_n\circ\alp}{1} \ + \
 \norm{\gb'_n\circ\alp}{1},
\label{threshold.convergence.e1}
\end{eqnarray}
where $\gs'_n :=  \lb|\gs-\gs_n\rb|$ and
$\gb'_n : =  \lb|\gb-\gb_n\rb|$.

For any $\del>0$, 
let $\bW_\del^s$ and $\bW_\del^b$ be as prior to Lemma \ref{conv.lemma.2}.
If
$\Del_n^s:=\CC{s_0,s_1}\symdif\CC{s_0^n,s_1^n}$,
then $\gs'_n  =  \chr{\Del_n^s}$.
If $n$ is big enough, then $|s^n_0-s_0|<\del$ and  $|s^n_1-s_1|<\del$,
so $\Del_n^s \subset \bW^s_\del$.  Thus,
\beqn
\norm{\gs'_n\circ\alp}{1}
 \quad =  \quad \lam\lb[\alp^{-1}(\Del_n^s)\rb]
 \quad \leq  \quad \lam\lb[\alp^{-1}(\bW_\del^s)\rb]
\quad=\quad M_\ba^s(\del).
\label{threshold.convergence.e2}
\eeqn
Likewise, if $\Del_n^b:=\CC{b_0,b_1}\symdif\CC{b_0^n,b_1^n}$,
then $\gb'_n  =  \chr{\Del_n^b}$.  
If $n$ is big enough, then
$|b^n_0-b_0|<\del$ and  $|b^n_1-b_1|<\del$,
so $\Del_n^b \subset \bW^b_\del$.  Thus,
\beqn
\norm{\gb'_n\circ\alp}{1}
 \quad =  \quad \lam\lb[\alp^{-1}(\Del_n^b)\rb]
 \quad \leq  \quad \lam\lb[\alp^{-1}(\bW_\del^b)\rb]
\quad=\quad M_\ba^b(\del).
\label{threshold.convergence.e3}
\eeqn
Combining equations 
(\ref{threshold.convergence.e1}) to (\ref{threshold.convergence.e3}), we get
\[
\norm{\Life_n(\ba)-\Life(\ba)}{1}
\quad\leq\quad
M_\ba^s(\del) \ + \ M_\ba^b(\del)
\quad\goesto{\del\goto0}{}\quad 0,
\]
because $M(\ba)=0$, because $\ba\in\oAX$ by hypothesis.

{\bf(b)}\quad  The set $\oAX$ is $\shift{}$-invariant.
Given part {\bf(a)} together with
Proposition \ref{fixed.point.convergence.lemma}, it suffices to
show that $\{\Life_n\}_{n=1}^\oo$ is eventually $\Ll$-equicontinuous at
every $\ba\in\oAX$.

  To do this, fix $\ba\in\oAX$, and $\eps>0$.  We want  $\del>0$
and $N\in\Natur$ so that, if $\ba'\in\lAX$ and $\norm{\ba-\ba'}{1}<\del$, then
$\norm{\Life_n(\ba)-\Life(\ba')}{1}<\eps$ for all $n\geq N$.

Let $K:=\norm{\gk}{\oo}$. 
Define measure $\mu_\ba$ on $\Real$ by
$\mu_\ba[\bU] :=
 \lam\lb[\alp^{-1}(\bU)\rb]$ for any measurable subset $\bU\subset\Real$.
 For any $n\in\Natur$, and any $\del>0$,
let ${}^n M^s_\ba(\del) \ :=\ \mu_\ba\lb[{}^n\bW^s_\del \rb]$,
where
${}^n\bW^s_\del  :=  \OO{s^n_0-\del,s^n_0+\del}\union
\OO{s^n_1-\del,s^n_1+\del}$, and
let
${}^nM^b_\ba(\del)  :=   \mu_\ba\lb[{}^n\bW^b_\del\rb]$,
where
${}^n\bW^b_\del  :=  \OO{b^n_0-\del,b^n_0+\del}\union
\OO{b^n_1-\del,b^n_1+\del}$.
If $\norm{\ba-\ba'}{1}<\del$, 
then Lemma  \ref{conv.lemma.2}(a,b,c) says
\beqn
\label{threshold.convergence.e4}
\hspace{-2em}
\mbox{$\forall \ n\in\Natur$,} \qquad
\norm{\Life_n(\ba)-\Life_n(\ba')}{1} \quad\leq\quad
2\del \ + \  {}^n M^s_\ba(K\del) \ + \ {}^n M^b_\ba(K\del).
\qquad\eeqn
For any $\del>0$, let 
$M_\ba^s(\del):=\mu_\ba[\bW^s_\del]$ and $M_\ba^b(\del):=\mu_\ba[\bW^b_\del]$.
Lemma  \ref{conv.lemma.2}(d) says $\D\lim_{\del\goto0} M_\ba^s(\del)
 =  0$, so if $\del$ is small enough, then
 $M_\ba^s(2K\del)<\eps/4$.
Find $N_0\in\Natur$ so that, for all $n \geq N_0$,
 \  $|s_0^n - s_0|< K\del$ and $|s_1^n - s_1|< K\del$.
  Then
${}^n\bW_{K\del}^s \subset \bW^s_{2K\del}$, so:
\beqn
\label{threshold.convergence.e5}
\hspace{-2em}
\forall \ n\geq N_0, \qquad
{}^n M^s_\ba(K\del)\  = 
\  \mu_\ba[{}^n\bW_{K\del}^s] \  \leq \  \mu_\ba[\bW_{2K\del}^s]
\ = \ M_\ba^s(2K\del) \ < \ \eps/4.\qquad
\eeqn
Likewise, Lemma  \ref{conv.lemma.2}(d) says that  if $\del$ is small enough, then
$M_\ba^b(2K\del)<\eps/4$. If $N_1$ is big enough, 
 then for all $n\geq N_1$,
 \  $|b_0^n - b_0|< K\del$ and $|b_1^n - b_1|< K\del$,
so
${}^n\bW_{K\del}^b \subset \bW^b_{2K\del}$, so:
\beqn
\label{threshold.convergence.e6}
\hspace{-2em}
\forall \ n\geq N_1, \qquad
{}^n M^b_\ba(K\del) \  =\ 
 \mu_\ba[{}^n\bW_{K\del}^b] \  \leq \ 
 \mu_\ba[\bW_{2K\del}^b]
\ = \ M_\ba^b(2K\del) \  < \ \eps/4.
\qquad
\eeqn
Let $N:=\max\{N_0,N_1\}$.
Thus, combining equations (\ref{threshold.convergence.e4})
to (\ref{threshold.convergence.e6}), we get:
\[
\hspace{-2em}
\mbox{For all $n\geq N$,} \qquad
\norm{\Life_n(\ba)-\Life_n(\ba')}{1}
\quad<\quad
2\del + \eps/4 + \eps/4 \quad=\quad  2\del + \eps/2
\]
Assume $\del <\eps/4$. If $\norm{\ba-\ba'}{1}<\del$,
then 
$\norm{\Life_n(\ba)-\Life_n(\ba')}{1} <\eps$ for all $n \geq N$.
\ethmprf

\section{From {\em Larger Than Life} to {\em RealLife}
\label{S:LtL}}

  In  what sense is a {\em RealLife} EA the `continuum limit' of a
sequence of {\em Larger than Life} CA?  We'll
construct a `discrete approximation' $\eLife$ of $\Life$ (for any
$\eps>0$), which is isomorphic to an LtL CA
with radius of order $\sO(1/\eps)$ (Proposition \ref{LtL.2.Reallife}).
We will then prove:

\Theorem{\label{pointwise.convergence}}
{
Fix $(s_0,b_0,s_1,b_1)\in\Theta$ and $\gk\in\sK$.
Let $\Life$ be the resulting RealLife EA.
\bthmlist
\item
  If $\ba\in\oAX$, then 
\ $\D \Lllim_{\eps\goto0}\, \eLife(\ba) \ = \ \Life(\ba)$.

\item  If $\D \lim_{n\goto\oo}\,\eps_n  =  0$, 
then $\{\eLife[\eps_n]\}_{n=1}^\oo$ evolves to
$\Life$ on $\oAX$.
\ethmlist
}

{\em Remarks:} (a)  It is clearly impossible to exactly simulate
a {\em RealLife} EA on a digital computer; the best we can do is simulate a
large-radius {\em Larger than Life} CA.
Theorem \ref{pointwise.convergence}(a) guarantees this is
will yield a `good approximation' of {\em RealLife}. 

(b) Evans \cite{Evans3} has found that LtL CA of
increasingly large radii have life forms
which are virtually identical after rescaling [see Figure
\ref{fig:four.bugs}].  This, combined with Theorem
\ref{pointwise.convergence}(b), provides compelling evidence (but
not proof) that {\em RealLife} EA have life forms 
which are morphologically similar to those seen in large-scale
LtL CA.\hfill$\diamondsuit$

\breath

To start, fix $\eps>0$, and let $\eZD \ := \
\set{\eps\fz}{\fz\in\ZD}$.  Let $\sB_\eps$ be the sigma algebra
generated by all $D$-dimensional half-open cubes of sidelength $\eps$,
with centres in $\eZD$.  That is, $\sB_\eps$ is generated by
$\set{\bC(\fz,\eps)}{\fz\in\ZD}$, where, for any $\fz\in\ZD$,
$\bC(\fz,\eps) := \eps\fz+\CO{\frac{-\eps}{2},\frac{\eps}{2}}^D$ is
the half-open cube of sidelength $\eps$, centred at $\eps\fz$.  If
$\eLl:=\Ll(\dX,\sB_\eps,\lam)$, then $\eLl\subset\Ll\intsct\Loo$.  Let
$\eAX := \lAX\intsct\eLl$.  Let $\elll:=\ell^1(\ZD)$.  For any
$\ba\in\eLl$, we define $\tlba\in\ell^1$ by $\tla_\fz := \ba(\eps\fz)$
for all $\fz\in\ZD$.  Thus, $\ba(c)=\tla_\fz$, for all $c\in \bC(\fz,\eps)$.
It follows:

\Lemma{\label{banach.iso1}}
{
For any $\ba\in\eLl$, 

{\bf(a)}
  $\norm{\ba}{\oo} \ = \ \norm{\tlba}{\oo}$, \quad
and \quad
{\bf(b)}  $\norm{\ba}{1} \ = \ \eps^D\cdot\norm{\tlba}{1}$. 

  Thus,
\bthmlist \setcounter{enumi}{2}
\item The map $\eLl \ni\ba\mapsto\tlba\in  \elll$
is a Banach space isomorphism.

\item Let $\lAZD := \AZD\intsct\elll$ be the set of $\ZD$-indexed
configurations with finite support.  Then
$\eAX \ni \ba\mapsto \tlba\in \lAZD$ is bijection.\qed
\ethmlist
}

  Let $\sM(\eZD)$ be the Banach space of signed measures on $\eZD$. 
For any $\fz\in\ZD$, let $\del_{\eps\fz}$ be the point mass at $\eps\fz$.
If $\gk\in\sM(\eZD)$, then $\D  \gk \ = \ 
 \sum_{\fz\in\ZD} k_\fz \del_{\eps\fz}$,
 for some $\tlgk:=[k_\fz]_{\fz\in\ZD}\in\elll$.
This defines a Banach isomorphism $\sM(\eZD)\ni\gk\mapsto\tlgk\in \elll$.
It is easy to show:

\Lemma{\label{banach.iso2}}
{
  If $\gk\in\sM(\eZD)$ and $\ba\in\eLl$, then: 
\quad {\bf(a)}  \  $\gk*\ba \in \eLl$, \quad
and \\ \quad {\bf(b)}  $\widetilde{\gk*\ba}  =  \tlgk * \tlba$. \qed
}

  If $\gk\in\sK$, we define $\ebgk\in\sM(\eZD)$ by
\beqn
\label{ebgk}
\hspace{-2em}
  \ebgk \ \ :=\ \  \sum_{\fz\in\ZD} k_\fz \del_{\eps\fz},
\ \ \mbox{where} \ \
k_\fz \ := \  \int_{\bC(\fz,\eps)} \gk(c) \ d\lam[c],
\ \ \mbox{for any $\fz\in\ZD$.}
\quad
\eeqn
We then define $\etgk := \widetilde{\ebgk} \ = \ 
[k_\fz]_{\fz\in\ZD} \in\elll\subset\ell^\oo(\ZD)$.
  It is easy to check that
\beqn
\label{etgk.vs.gk}
\begin{array}{rrcl}
{\bf(A)} &
\norm{\etgk}{1} & = &\D \sum_{\fz\in\ZD} k_\fz \ = \ \int_\dX \gk \ = \ 1, \\
\qquad\And\quad
{\bf(B)} & \norm{\etgk}{\oo}& \leq&  \eps^D \cdot \norm{\gk}{\oo}.
\end{array}
\eeqn
  If $\ba$ is $\sB_\eps$-measurable, then Lemma 
\ref{banach.iso2}(a) says $\alp:=\ebgk*\ba$ is also $\sB_\eps$-measurable, so
$\gs\circ\alp$ and $\gb\circ\alp$ are also $\sB_\eps$-measurable.  Thus,
we can define $\ebLife:\eLl\into\eLl$ by
\beqn
\label{elife.defn}
 \ebLife(\ba)\quad:=\quad
 \ba\cdot \gs\circ(\ebgk * \ba) \ + \  (1-\ba)\cdot \gb\circ(\ebgk * \ba),
\qquad\forall \ba\in\eLl.
\eeqn
 Note that $\ebLife(\eAX)\subseteq \eAX$
so $\ebLife$ restricts to a transformation on $\eAX$.

  Let $\dE_\eps:\Ll\into\eLl$ be the conditional expectation operator
for the sigma-algebra $\sB_\eps$.  That is:
\[
\mbox{For any $\ba\in\Ll$ and $x\in\dX$,}\qquad
  \dE_\eps[\ba](x) \quad:=\quad \frac{1}{\eps^D} \int_{\bC} \ba(c) \ d\lam[c],
\]
where $\bC\in\sB_\eps$ is the unique $\eps$-cube containing $x$. 
   We extend $\ebLife$ to
a function $\eLife:\Ll\into\eLl$ by defining $\eLife(\ba)
:= \ebLife(\dE_\eps[\ba])$ for any $\ba\in\Ll$.  Note that
$\eLife(\ba)=\ebLife(\ba)$ for any $\ba\in\eLl$, because
$\dE_\eps$ acts as the identity on $\eLl$.   Thus, we 
will suppress the distinction between $\eLife$ and $\ebLife$, and
write both as ``$\eLife$''.
  
\breath

If $\etgk := \widetilde{\ebgk}\in\elll$, then
we  define the operator $\etLife:\elll\into\elll$ by
\[
 \etLife(\ba)\quad:=\quad
 \ba\cdot \gs\circ(\etgk * \ba) \ + \  (1-\ba)\cdot \gb\circ(\etgk * \ba),
\qquad\forall \ba\in\elll.
\]

\Proposition{\label{LtL.2.Reallife}}
{
\bthmlist
\item $\etLife(\AZD)\subseteq\AZD$, and $\etLife:\AZD\into\AZD$ is
a Larger than Life CA, with birth interval $\CC{b_0,b_1}$ and
survival interval $\CC{s_0,s_1}$.

 In particular, suppose $\gk=\lam[\dK]^{-1}\chr{\dK}$, where
$\dK=\CC{-1,1}^D\subset\dX$.  If $\eps=1/n$, let 
$\dK_n := \CC{-n...n}^D\subset\ZD$.  Then $\etLife$ is
an LtL CA with neighbourhood $\dK_n$.

\item The map $\eLl\ni\ba\mapsto \tlba\in\elll$ is
a dynamical isomorphism 
$(\eLl,\eLife)\cong(\elll,\etLife)$, which restricts to 
a dynamical isomorphism  $(\eAX,\eLife)\cong(\lAZD,\etLife)$.
\ethmlist
}
\bthmprf {\bf(a):} Use eqn.(\ref{etgk.vs.gk}A).  
 {\bf(b):} \  Lemma \ref{banach.iso2}(b) implies 
 $\widetilde{\eLife(\ba)} \ = \ \etLife(\tlba), \ 
\forall \ \ba\in\eLl$.
\ethmprf

{\em Remarks:}
 By choosing $s_0\leq b_0 < b_1 \leq s_1$ and $\gk$ appropriately,
 we can obtain any LtL CA in part (a) of Proposition \ref{LtL.2.Reallife}.
Part (b) means that
 each life form of $(\eAX,\eLife)$
corresponds to a life form of 
the LtL CA $(\AZD,\etLife)$.\hfill$\diamondsuit$

\breath

\ignore{
\Lemma{\label{conv.lemma.0}}
{
Fix $\ba\in\bAX$.  If
 $\eps$ is small enough, and $\eba:=\dE_\eps[\ba]$, then
$\norm{\eba - \ba}{1} \ \leq \  2L(\ba)\cdot\eps$.
}
\bthmprf
Suppose $\ba=\chr{\bA}$ for some $\bA\subset\dX$, and let $L:=L(\ba)$.  Let
$\bA_\eps\in\sB_\eps$ be the largest $\sB_\eps$-measurable set
contained in $\bA$ and let $\compl{\bA}_\eps\in\sB_\eps$ be the largest
$\sB_\eps$-measurable set contained in $\compl{\bA}$.  Let
$\partial_\eps \bA := \dX\setminus (\bA_\eps \disj
\compl{\bA}_\eps)$.  Thus,
\[
\norm{\eba - \ba}{1} \quad=\quad
 \int_{\dX} |\ba(x) - \eba(x)| \ d\lam[x]
 \quad \eeequals{(*)} \quad
 \int_{\partial_\eps \bA} |\ba(x) - \eba(x)| \ d\lam[x]
 \quad \leeeq{(\dagger)} \quad \lam[\partial_\eps \bA] 
\quad\leeeq{(\ddagger)}\quad 2L\eps.
\]
$(*)$ If $x\in\bA_\eps$ or $x\in\compl{\bA}_\eps$,
 then $\ba(x) \ = \ \eba(x)$.
\quad 
$(\dagger)$ If $x\in\partial_\eps\bA$, then $|\ba(x) - \eba(x)|\leq 1$.

$(\ddagger)$  $\partial_\eps \bA$ is  a covering of $\partial \bA$
by $\eps$-cubes from $\sB_\eps$.   This covering needs approximately
$L/\eps^{D-1}$ cubes, each of volume $\eps^D$.  Thus,
$\D\lim_{\eps\goto0} \frac{1}{\eps} \lam[\partial_\eps \bA]
 =  L$, so if $\eps$ is small, then $\lam[\partial_\eps \bA] < 2L\eps$.
\ethmprf
}

Let 
  $\bAX \ := \ 
 \set{\ba\in\oAX}{\mbox{$\ba=\chr{\bA}$ for some compact
$\bA\subset\dX$ with $\lam[\partial\bA] = 0$}}$.

\Lemma{\label{bAX.dense.in.oAX}}
{
$\bAX$ is a $\bL^1$-dense subset of $\oAX$.
}
\bthmprf
  Let $\ba\in\oAX$, and suppose $\ba=\chr{\bA}$ for some measurable
$\bA\subset\dX$.  For any $\eps>0$, there is a finite disjoint union
$\bU\subset\dX$ of bounded closed cubes, such that $\lam[\bA\symdif\bU]<\eps$ 
\cite[Thm. 2.40(c), p.68]{Folland}.
If $\bu=\chr{\bU}$, then
$\norm{\ba-\bu}{1} = \lam[\bA\symdif\bU]<\eps$.
Let $\upsilon := \gk*\bu$; by slightly increasing/decreasing the
sidelengths of the cubes of $\bU$, we can slightly increase/decrease
$\ups$ everywhere, and thereby ensure that $M(\bu):=
\lam\lb[\ups^{-1}\{s_0,b_0,b_1,s_1\}\rb]=0$
(the argument is similar to the `density' proof
in Theorem \ref{oAX.is.comeager}).
  Hence, $\bu\in\oAX$.  Also, $\bu\in\bAX$, because
$\lam[\partial\bU]=0$, because
$\partial\bU$ is a finite union of $(D-1)$-dimensional cube faces.
\ethmprf

\Lemma{\label{conv.lemma.1}}
{
  Let $\ba\in\lAX$ and  $\gk\in\sK$.  For any $\eps>0$, let
$\eba:=\dE_\eps[\ba]$.  Let $\ebgk\in\sM(\eZD)$ be as
in {\rm(\ref{ebgk})}.  Then
\bthmlist
  \item  $\D \lim_{\eps\goto0} \norm{\eba-\ba}{1} \ = \ 0$.

 \item If $\ba\in\bAX$,
then  $\D\lim_{\eps\goto 0} \norm{\ebgk*\eba - \gk*\ba}{\oo} \ = \ 0$.
\ethmlist
}
\bthmprf {\bf(a)} is basic martingale theory (see \cite[Thm.3.3.2]{Borkar} or
\cite[Cor.5.2.7]{Strook}).  

 {\bf(b)} \ First observe that
\beqn
\hspace{-2em}\norm{\ebgk*\eba - \gk*\ba}{\oo} 
\quad\leq\quad
\norm{\ebgk*\eba - \gk*\eba}{\oo}
\ + \ 
  \norm{\gk*\eba - \gk*\ba}{\oo}.
\qquad\label{conv0}
\eeqn
  Let  $K:=\norm{\gk}{\oo}$.  Then
$\norm{\gk*\eba - \gk*\ba}{\oo}
\ = \
 \norm{\gk*(\eba - \ba)}{\oo}
\ \leq \ K\cdot\norm{\eba - \ba}{1}$,
by Young's inequality.  Thus
{\bf(a)} implies that $\D\lim_{\eps\goto0} \, \norm{\gk*\eba - \gk*\ba}{\oo}
 \ = \ 0$.
Thus, \eref{conv0} implies that
\[
\lim_{\eps\goto 0} \, \norm{\ebgk*\eba - \gk*\ba}{\oo}
\quad\leq\quad
\lim_{\eps\goto 0} \, \norm{\ebgk*\eba - \gk*\eba}{\oo}.
\]  
\newcommand{\distal}[1]{\stackrel{#1}{\leadsto}}
We claim the right-hand limit is zero.
To see this, we need some notation.
For any $n\in\Natur$, and $x,y\in\dX$, we write ``$x \distal{n} y$''
to mean that $y \in x + \CO{-\frac{n\eps}{2},\frac{n\eps}{2}}^D$ (note
that this relation is not quite symmetric).
If $\bU,\bV\subset\dX$, then ``$\bU\distal{n}\bV$'' means
$u\distal{n}v$ for some $u\in\bU$ and $v\in\bV$.
For example, if $\fz\in\ZD$, recall that $\bC(\fz,\eps) := 
\eps\fz+\CO{-\frac{\eps}{2},\frac{\eps}{2}}^D
\ = \  \set{c\in\dX}{\eps\fz \distal{1} c}$.

{\em Note:}  if $x\not\distal{2}y$, then $x$ and $y$ cannot be in the same
$\sB_\eps$-cube 
[because if $\bB\in\sB_\eps$ and $x,y\in\bB$, then $x\distal{2}y$].

Suppose $\ba=\chr{\bA}$, for some $\bA\subset\dX$.
For any $x\in\dX$, let
$\bW_\eps^x :=\set{\fz\in\ZD}{ \partial\bA \distal{3} (x-\eps\fz)}$
and  let $\bV_\eps:=\ZD\setminus\bW^x_\eps$. Then
\begin{eqnarray}
\nonumber
  \ebgk*\eba(x)
&=&
\sum_{\fz\in\ZD} k_\fz \cdot \eba(x-\eps\fz)
\\ \hspace{-2em} &=&
\sum_{\fw\in\bW^x_\eps} k_\fw \cdot \eba(x-\eps\fw)
\ + \
\sum_{\fv\in\bV_\eps} k_\fv \cdot \eba(x-\eps\fv),\qquad
\label{conv1}
\end{eqnarray}
and
\begin{eqnarray}
\nonumber 
\hspace{-2em}
\lefteqn{  \gk*\eba(x)
\quad=\quad
\int_{\dX} \gk(y)\cdot \eba(x-y) \ dy
\quad=\quad 
\sum_{\fz\in\ZD}  
\int_{\bC(\fz,\eps)} \gk(c)\cdot \eba(x-c) \ dc} \\
\hspace{-2em}
&=&
\sum_{\fw\in\bW^x_\eps}  \int_{\bC(\fw,\eps)} \gk(c)\cdot \eba(x- c) \ dc
\ + \ 
\sum_{\fv\in\bV_\eps} 
 \int_{\bC(\fv,\eps)} \gk(c)\cdot \eba(x- c) \ dc.
\qquad\qquad
\label{conv2}
\end{eqnarray}

\Claim{\label{conv3}
For any $\fv\in\bV_\eps$, \
 $\D \int_{\bC(\fv,\eps)} \gk(c)\cdot \eba(x- c) \ dc 
\ = \ \eba(x-\eps\fv)\cdot k_\fv$.}
\bclaimprf
Let $\bC:=x-\bC(\fv,\eps)$.
If $\fv\in\bV_\eps$, then $\partial\bA\not\distal{3} (x-\eps\fv)$.
Thus, $\partial\bA\not\distal{2} \bC$.  Thus if $\bB\in\sB_\eps$
and $\bB\intsct\bC\neq\emptyset$, then $\bB\intsct\partial\bA=0$
[because if $b\in\bB$, then $b \distal{2} \bC$, so $b\not\in\partial\bA$].
Thus, no $\sB_\eps$-cube can intersect both $\bC$ and
$\partial\bA$.  
Thus, $\eba$ is constant on $\bC$, because:
\bitem
\item[{\em either}] ${\eba}\restr{\bC} \equiv 0$, if every element of $\bC$ lies in
a  $\sB_\eps$-cube entirely {\em outside} $\bA$, 
\item[{\em or}] ${\eba}\restr{\bC}\equiv 1$,
if every element of $\bC$ lies in a $\sB_\eps$-cube
that is  entirely {\em inside} $\bA$.
\eitem
\beq
\hspace{-2em}
\mbox{Thus},\quad \int_{\bC(\fv,\eps)} \gk(c)\cdot \eba(x-c) \ dc
\ &\eeequals{(\dagger)}& \ 
 \int_{\bC} \gk(x-c')\cdot \eba(c') \ dc'
\\ &=& \ \ 
 \eba(x-\eps\fv) \int_{\bC(\fv,\eps)}  \gk( c)  \ dc
\ \ \eeequals{(*)} \ \ 
 \eba(x-\eps\fv) \cdot k_\fv, 
\eeq
where $(\dagger)$ is the change of variables $c':= x-c$, and
 $(*)$ is by eqn.(\ref{ebgk}).
\eclaimprf
Applying Claim \ref{conv3} to eqn.(\ref{conv2}), we conclude that
\beqn
\hspace{-1em}
 \gk*\eba(x)
\quad=\quad
\sum_{\fw\in\bW^x_\eps}
      \int_{\bC(\fw,\eps)} \gk(c)\cdot \eba(x- c) \ dc
\ + \ 
\sum_{\fv\in\bV_\eps}  k_\fv\cdot \eba(x-\eps\fv).
\qquad\label{conv4}
\eeqn
Subtracting eqn.(\ref{conv4}) from eqn.(\ref{conv1}), we see that
\begin{eqnarray}
\nonumber
\lefteqn{\hspace{-1.5em} \lb|\ebgk*\eba(x) - \gk*\eba(x) \rb|
\quad\leq\quad
\sum_{\fw\in\bW^x_\eps} \lb|k_\fw \cdot \eba(x-\eps\fw) -  
\int_{\bC(\fw,\eps)} \gk(c) \cdot \eba(x-c) \ dc\rb|} \\
&\leq& \ \ 
\sum_{\fw\in\bW^x_\eps} \lb(\maketall \norm{\etgk}{\oo}\cdot \norm{\eba}{\oo}
+   \lam[\bC(\eps)]\cdot \norm{\gk}{\oo}\cdot\norm{\eba}{\oo} \rb)
\nonumber \\&\leeeq{(*)}& \ \ 
2K\cdot\#(\bW^x_\eps)\cdot\eps^D.\qquad\qquad
\label{conv5}
\end{eqnarray}
$(*)$ is by eqn.(\ref{etgk.vs.gk}B),
the fact that $\lam[\bC(\eps)] = \eps^D$, and the fact that
$\norm{\eba}{\oo} \leq 1$, because $\eba = \dE_\eps[\chr{\bA}]$.
It remains to control the right hand side of eqn.(\ref{conv5}).

\Claim{\label{conv7}
Fix $\eps>0$.   For any $x,y\in\dX$, \ $\#(\bW^x_\eps) =  \#(\bW^y_\eps)$.}
\bclaimprf
Fix $x\in\dX$, and let $\bA_x := x-\bA = \set{x-a}{a\in\bA}$.
Then $\partial\bA_x := x-\partial\bA$.
Thus,
\beqn
\label{conv8}
\bW^x_\eps\quad=\quad
\set{\fz\in\ZD}{\eps\fz\distal{3}\partial\bA_x},
\eeqn
because for any $\fz\in\ZD$,
\[
\statement{$\fz\in\bW^x_\eps$}
 \ \iff \ 
\statement{$\partial\bA\distal{3}(x-\eps\fz)$}
\ \iff \
\statement{$\eps\fz\distal{3}(x-\partial\bA)=\partial\bA_x$}.
\]
Define $\xi:\Real\into\dX$ by
$\xi(t):=x+(\eps t,0,...,0) $;  hence $\xi(0)=x$.  Let 
$\bW_t:=\bW^{\xi(t)}_\eps  =  
\set{\fz\in\ZD}{\eps\fz\distal{3}\partial\bA_{\xi(t)}}$.
Hence $\bW_0=\bW^x_\eps$.
As we increase $t$, elements of $\Zahl^D$ enter $\bW_t$ at the same
rate as they leave.  To be precise, for each
 $\fz\in\Zahl^D$, let 
\[
T_i(\fz) \ := \ \inf\set{t\in\Real}{\fz\in\bW_t}
\ \ \And \ 
T_o(\fz)\ := \ \max\set{t\in\Real}{\fz\in\bW_t}
\]
be the `entrance time' of $\fz$ into $\bW_t$, and `exit time' 
out of $\bW_t$, respectively. Also, for any  $t >0$, let
\[\hspace{-2em}
  \bI(t) \ := \ \set{\fz\in\Zahl^D}{0\leq T_i(\fz)\leq t}
\ \And \ 
  \bO(t) \ := \ \set{\fz\in\Zahl^D}{0\leq T_o(\fz)\leq t}.
\]
Thus, $\#(\bW_t) \ = \ \#(\bW_0) + \#\bI(t) - \#\bO(t)$.

But for all $t>0$, \ $\#\bI(t)=\#\bO(t)$.  To see this, let
$\fz\in\Zahl^D$;  if $\fz' = \fz+(3,0,\ldots,0)$, then clearly $T_o(\fz) = T_i(\fz')$. This yields a bijection $\bI(t)\ni\fz\mapsto\fz'\in\bO(t)$.  

  It follows that $\#(\bW_t)  =  \#(\bW_0)$ for all $t>0$.
In other words, $\#\bW^x_\eps$ is constant as we vary the first
coordinate of $x$.  The same applies to any other coordinate.
\eclaimprf

\Claim{\label{conv9}
For any $\del>0$, there is  $\eps>0$ such that,
for all $x\in\dX$, \ $\#(\bW^x_\eps)\cdot\eps^D < \del/2K$.}
\bclaimprf
Fix $x\in\dX$ and let $\bA_x$ be as in Claim \ref{conv7}.  Then
eqn.(\ref{conv8}) says that
 $\bW^x_\eps$ is the
set of all cube centres of a covering $\bU_\eps$
of $\partial\bA_x$ by (overlapping) $\frac{3}{2}\eps$-cubes, defined
\[
  \bU_\eps \quad := \quad \Union_{\fw\in\bW^x_\eps} \bC'(\fw,\mbox{$\frac{3}{2}\eps$}),
\qquad
\mbox{where \ $\bC'(\fw,\frac{3}{2}\eps) \ := \
 \eps\fw + \CO{\frac{-3}{2}\eps, \frac{3}{2}\eps}^D$.}
\]
 Thus,
$\lam[\bU_\eps]  \ \stackrel{\scriptscriptstyle(*)}{\geq} \
\frac{1}{3^D}  \#(\bW^x_\eps)\cdot \lam[\bC'(0,\frac{3}{2}\eps)] 
\ = \ \#(\bW^x_\eps)\cdot(3\eps)^D/3^D \ = \  \#(\bW^x_\eps)\cdot \eps^D$,
where $(*)$ is because each point of $\bU_\eps$ is covered by at most $3^D$ distinct
$\frac{3}{2}\eps$-cubes in the covering.
 If $\bB_\eps:=\partial\bA+\OO{-3\eps,3\eps}^D$,
 then $\partial \bA_x\subset\bU_\eps\subset\bB_\eps$.  Thus,
\begin{eqnarray}
\lim_{\eps\goto0}\, \#(\bW^x_\eps)\cdot\eps^D  & \ \leq \ & 
\lim_{\eps\goto0}\, \lam[\bU_\eps]
\ \ \leq\ \ 
\lim_{\eps\goto0}\, \lam[\bB_\eps] 
\nonumber \\ &\ \eeequals{(*)}\ &
 \ \  \lam[\partial\bA_x] \ \  =\ \  \lam[\partial\bA] \quad  
\eeequals{(\dagger)}\quad  0.
\label{conv6}
\end{eqnarray}
Here, $(\dagger)$ is  because $\ba\in\bAX$,
and $(*)$ is by `continuity from above' of the measure $\lam$,
because 
$\{\bB_\eps\}_{\eps>0}$ is a decreasing family of open sets 
with $\Intsct_{\eps>0} \bB_\eps = \partial\bA$
(because: if $x\in \Intsct_{\eps>0} \bB_\eps$, then
$x$ is a cluster point of $\partial\bA$, which is a closed set,
 so $x\in\partial\bA$).

 For any $\del>0$, \eref{conv6} yields  $\eps>0$ such that
 $\#(\bW^x_\eps)\cdot \eps^D \ < \  \del/2K$.  But then
Claim \ref{conv7} implies 
$\#(\bW^y_\eps)\cdot \eps^D \ < \  \del/2K$ for any $y\in\dX$.
\eclaimprf
If $\del>0$ and $\eps>0$ are as in Claim \ref{conv9}, then for any $x\in\dX$, 
eqn.(\ref{conv5}) says that
\[
\lb|\ebgk*\eba(x) - \gk*\eba(x) \rb|
\quad\leq \quad 
2K\#(\bW^x_\eps)\cdot\eps^D 
\quad<\quad
\del.
\]
Thus,
$\norm{\ebgk*\eba - \gk*\eba}{\oo}
\ < \ \del$.  Let $\del\goto0$ to get
$\D \lim_{\eps\goto0}   \norm{\ebgk*\eba - \gk*\eba}{\oo} \ = \ 0$.
\ethmprf

\bthmprf[Proof of Theorem \ref{pointwise.convergence}:] 
\setcounter{claimcount}{0}
\Claim{ If $\ba\in\bAX$, then 
\ $\D \Lllim_{\eps\goto0}\, \eLife(\ba) \ = \ \Life(\ba)$.}
\bclaimprf
Let $\alp:=\gk*\ba$.  For any $\eps>0$,
 let $\baralp_\eps := \ebgk*\eba$, where $\eba=\dE_\eps(\ba)$.  Then
\begin{eqnarray}
\nonumber
\hspace{-3em}
\eLife(\ba) & = & \eba\cdot \gs\circ\baralp_\eps \ + \ (1-\eba)\cdot \gb\circ\baralp_\eps,
\quad\mbox{by \eref{elife.defn},} \\ \hspace{-3em}
\nonumber
\And\quad 
\Life(\ba) & = & \ba\cdot \gs\circ\alp \ + \ (1-\ba)\cdot \gb\circ\alp,
\qquad \ \mbox{by \eref{afterlife.CA2}.}\\ \hspace{-3em}
\nonumber
\mbox{Thus,} \
\eLife(\ba)-\Life(\ba) &=&
(\eba-\ba)\cdot \gs\circ\baralp_\eps  \ + \
\ba\cdot \lb(\gs\circ\baralp_\eps - \gs\circ\alp\rb) \\ \hspace{-3em}
\nonumber
&&\qquad \ + \
(\ba-\eba)\cdot \gb\circ\baralp_\eps  \ + \
 (1-\ba)\cdot \lb(\gb\circ\baralp_\eps - \gb\circ\alp\rb). \\ \hspace{-3em}
\nonumber
\mbox{So,} \ 
\norm{\eLife(\ba)-\Life(\ba)}{1} &\leq&
\norm{\eba-\ba}{1}\cdot \norm{\gs\circ\baralp_\eps}{\oo}  \ + \
\norm{\ba}{\oo}\cdot \norm{\gs\circ\baralp_\eps - \gs\circ\alp}{1}
\\ \hspace{-3em} \nonumber
&& \ \ 
\ + \ \norm{\ba-\eba}{1}\cdot \norm{\gb\circ\baralp_\eps}{\oo} 
 \ + \
 \norm{1-\ba}{\oo}\cdot \norm{\gb\circ\baralp_\eps - \gb\circ\alp}{1}. \\ \hspace{-3em}
&\leq &
2\cdot\norm{\eba-\ba}{1} \ + \ 
\norm{\gs\circ\baralp_\eps - \gs\circ\alp}{1} \ + \
 \norm{\gb\circ\baralp_\eps - \gb\circ\alp}{1}.  \qquad\quad
\label{pointwise.convergence.eqn1}
\end{eqnarray}
Now,  Lemma \ref{conv.lemma.1}(a) says $\D \lim_{\eps\goto0} \norm{\eba-\ba}{1} \ = \ 0$.  
Also, $\ba\in\bAX$, so Lemma \ref{conv.lemma.1}(b) says 
 $\D \lim_{\eps\goto0}\norm{\baralp_\eps-\alp}{\oo} \ = \ 0$.  
Thus, Lemma \ref{conv.lemma.2}(b,d) says  
 $\D \lim_{\eps\goto0}\norm{\gs\circ\baralp_\eps-\gs\circ\alp}{\oo} \ = \ 0$ and
 $\D \lim_{\eps\goto0}\norm{\gb\circ\baralp_\eps-\gb\circ\alp}{\oo} \ = \ 0$,
because $\ba\in\oAX$.

\noindent Combine these facts with \eref{pointwise.convergence.eqn1} to get
 $\D \lim_{\eps\goto0} \norm{\eLife(\ba)-\Life(\ba)}{1} \ = \ 0$.
\eclaimprf

\Claim{The family $\{\eLife\}_{\eps>0}$ is $\Ll$-equicontinuous at
every $\ba\in\oAX$.}
\bclaimprf
Let  $\ba\in\oAX$ and
fix $\gam>0$. We want   $\del>0$ so that, 
if $\norm{\ba-\ba'}{1}<\del$, then, for all $\eps>0$, \ 
$\norm{\eLife(\ba)-\eLife(\ba')}{1}<\gam$. 
Now, $M(\ba)=0$, so Lemma \ref{conv.lemma.2}(d) yields 
$\del>0$ with
\beqn
\label{elife.equicont.e0}
\del \ < \ \frac{\gam}{4},
\qquad
M_\ba^s(K\del) \ < \  \frac{\gam}{4},
\quad\And\quad
M_\ba^b(K\del) \ < \  \frac{\gam}{4}.
\eeqn
(Here, $K:=\norm{\gk}{\oo}$.)  
  Suppose $\ba'\in\Ll$, with $\norm{\ba-\ba'}{1}<\del$.
Let $\eba:=\dE_\eps(\ba)$ and $\eba':=\dE_\eps(\ba')$; then
$\norm{\eba-\eba'}{1}\leq \norm{\ba-\ba'}{1}<\del$, 
because $\dE_\eps$ is a bounded linear
operator on $\bL^1$ with $\norm{\dE_\eps}{}=1$.   Let $\tlba$ and $\tlba'$ be
the corresponding elements of $\elll$; then 
$\norm{\tlba-\tlba'}{1} = \norm{\eba-\eba'}{1}/\eps^D \ < \del/\eps^D$,
by Lemma \ref{banach.iso1}(b).
Thus,  if $\alp=\ebgk*\eba$ and $\alp' = \ebgk*\eba'$, then 
\begin{eqnarray}\nonumber
\norm{\alp - \alp'}{\oo} 
&\eeequals{(*)}&
\norm{\widetilde{\ebgk*\eba} - \widetilde{\ebgk*\eba'}}{\oo} 
\quad\eeequals{(\dagger)}\quad
\norm{\etgk*\tlba - \etgk*\tlba'}{\oo} 
\\\nonumber
&\leeeq{(\diamond)}& \ 
\norm{\etgk}{\oo} \cdot \norm{\tlba-\tlba'}{1}
\ \ \leq  \ \ 
\norm{\etgk}{\oo} \cdot \frac{\del}{\eps^D}
\\ &\leeeq{(\ddagger)}& \ 
\norm{\gk}{\oo}\cdot\eps^D \cdot \frac{\del}{\eps^D}
\quad=\quad
K \del. \label{elife.equicont.e1}
\end{eqnarray}
 Here, $(*)$ is by Lemma \ref{banach.iso1}(a), \ 
 $(\dagger)$ is by Lemma \ref{banach.iso2}(b), \ $(\diamond)$ is by Young's inequality,
and $(\ddagger)$ is by eqn.(\ref{etgk.vs.gk}B).
  It is easy to prove the analog of
Lemma \ref{conv.lemma.2}(a) for $\eLife$.  Combined with
Lemma \ref{conv.lemma.2}(b) and \eref{elife.equicont.e1}, this yields: 
\[
\hspace{-2em}
\forall \ \eps>0,\qquad
\norm{\eLife(\ba)-\eLife(\ba')}{1}
\quad \leq \quad
   2\del \ + \  M_\ba^s(K\del) \ + \  M_\ba^b(K\del)
\quad\leeeq{(*)}\quad
  \gam.
\]
Here, $(*)$ is by \eref{elife.equicont.e0}.
This works for all $\eps>0$, so
$\{\eLife\}_{\eps>0}$ is equicontinuous at $\ba$.
\eclaimprf

{\bf(a)} \quad Let $\ba\in\oAX$, and fix $\gam>0$.  Theorem 
\ref{reallife.is.continuous} and Claim 2 together
yield $\del>0$ so that, for any $\ba'\in\oAX$, if
$\norm{\ba'-\ba}{1}<\del$, then
\beqn
\hspace{-2em}
\norm{\Life(\ba')-\Life(\ba)}{1} \ < \ \frac{\gam}{3},\quad
\mbox{and, $\forall \ \eps>0$,} \quad
\norm{\eLife(\ba)-\eLife(\ba')}{1} \ < \ \frac{\gam}{3}.
\qquad\label{foo}
\eeqn
By Lemma \ref{bAX.dense.in.oAX}, find some $\ba'\in\bAX$
with $\norm{\ba'-\ba}{1}<\del$.  Finally, Claim 1 yields 
some $\sE>0$ so that, if $0<\eps<\sE$, then 
 $\norm{\eLife(\ba')-\Life(\ba')}{1}<\frac{\gam}{3}$.  Thus,
 if $0<\eps<\sE$, then 
\beq
\hspace{-2em} \lefteqn{\norm{\eLife(\ba)-\Life(\ba)}{1}}
\\\hspace{-2em}
&\qquad\leq\quad&
\norm{\eLife(\ba)-\eLife(\ba')}{1}
+ \norm{\eLife(\ba')-\Life(\ba')}{1}
+ \norm{\Life(\ba')-\Life(\ba)}{1}
\quad \leeeq{(\ref{foo})} \quad
\gam.
\eeq
  Since this works for any $\gam$, we conclude that
$\D \Lllim_{\eps\goto0}\, \eLife(\ba) \ = \ \Life(\ba)$.

{\bf(b)}  follows from part {\bf(a)}, Claim 2, and
Proposition \ref{fixed.point.convergence.lemma}.
\ethmprf

\section{A gallery of still lifes \label{S:rot.fixed.point}}

\Proposition{\label{fixed.point.cond}}
{
Let $\ba=\chr{\bA}\in\AX$,
let $\gs\circ(\gk*\ba)=\chr{\bS}$ and let $\gb\circ(\gk*\ba)=\chr{\bB}$,
for some $\bA,\bS,\bB\subset\dX$.  Then
\quad
$ \statement{$\Life(\ba)=\ba$}
\iff
\statement{$\bB\subseteq\bA\subseteq\bS$}$.
}
\bthmprf
 $\Life(\ba)=\chr{\bU}$, where 
$\bU := (\bA\intsct\bS) \disj (\compl{\bA}\intsct\bB)$.  Thus,

\noindent$\begin{array}[b]{rcl}
\statement{$\Life(\ba)=\ba$}
&\iff&
\statement{$\bA = \bU = (\bA\intsct\bS) \disj (\compl{\bA}\intsct\bB)$}
\\&\iff&
\statement{$\bA\intsct\bS =\bA$ and $\compl{\bA}\intsct\bB = \emptyset$}
\ \iff\ 
\statement{$\bB\subseteq\bA\subseteq\bS$}.
\end{array}$
\ethmprf

\Corollary{\label{no.half.life}}
{
   If $\supp{\gk}$ is a neighbourhood of zero,
and $s_0\geq\frac{1}{2}$, then $\Life$ has no bounded still lifes.
}
\bthmprf
  Suppose $\bH\subset\dX$ was any open halfspace such that $0\in\partial\bH$.
Then
\[
\int_\bH \gk(-h) \ d\lam[h]
\ + \  \int_{\compl{\bH}} \gk(-h) \ d\lam[h] \quad = \quad 
\int_\dX \gk(-x) \ d\lam[x] 
\quad=\quad 1.
\]
Thus (replacing $\bH$ with $\compl{\bH}$ if necessary)
we can assume $\int_\bH \gk(-h) \ d\lam[h] \leq \frac{1}{2}$.

 Let $\bA\subset\dX$ be some bounded set, and let
 $\ba:=\chr{\bA}\in\AX$.  Let $\bC$ be the convex closure
of $\bA$.  By translating $\bA$ if necessary, we can assume that
$\bC\subset\bH$, and that $\partial\bC$ is tangent to $\partial\bH$.
By slightly rotating $\bH$ if necessary, we can assure
that $\partial\bC$ is tangent to $\partial\bH$ 
at precisely one extremal point $e$, while still
preserving the inequality $\int_\bH \gk(-h) \ d\lam[h] \leq \frac{1}{2}$.
If $e$ is extremal in $\bC$, then $e\in\bA$.
By further translating $\bA$, we  assume $e=0$.  Thus,
\beq
 \gk*\ba(0) &\ = \ &
\int_\bA \gk(-a) \ d\lam[a]
 \ \   \leeeq{(\ddagger)} \ \   \int_\bC \gk(-c) \ d\lam[c]
\ \   \lt{(*)} \ \  
\int_\bH \gk(-h) \ d\lam[h]
\\&  \ \leq \ &  
 \frac{1}{2}
\ \ \leq \ \  s_0. 
\eeq
Here, $(\ddagger)$ is because $\bA\subseteq\bC$.
To see $(*)$, let $\dK=\supp{\gk}$, a neighbourhod of $0$. But
 $0$ is a cluster point of $(\bH\setminus\bC)$
(because $\bC\intsct\partial\bH \ = \{0\}$),
so $\bU:=-\dK\intsct(\bH\setminus\bC)$ is a nonempty open set.
Thus
\[
\int_\bH \gk(-h) \ d\lam[h] \ - \ \int_\bC \gk(-c) \ d\lam[c]
\quad = \quad 
\int_{\bU} \gk(-u) \ d\lam[u] 
\quad>\quad 0.
\]
 Thus $\gs\circ(\gk*\ba)(0)=0$, so (in terms of
Proposition \ref{fixed.point.cond})
 $0\not\in\bS$.   Hence, $\bA\not\subseteq\bS$;
so $\chr{\bA}$ can't be a still life.
\ethmprf

  In a sense, as $s_0$ becomes larger, the maximum convex curvature of
the boundary of a still life becomes smaller.  For example, suppose
$D=2$.  If $s_0>\frac{1}{4}$, then no still life can have a convex
right angle or  acute angle on its boundary (because if $\ba$
was a still life whose boundary made an angle $\leq\pi/2$ at $x$, then
$\gk*\ba(x)\leq\frac{1}{4}<s_0$).  If $s_0>\frac{1}{2}$, then the
boundary of a still life must be concave everywhere, which is
impossible; hence Corollary \ref{no.half.life}.

\breath

 Evans \cite{Evans2,Evans3} has found compact, `ball'-shaped still
lifes in many LtL CA, reminiscent of the well-known $2\x2$ square {\em
block} from Conway's {\em Life}.  We'll now construct a broad family
of such still lifes.  If $\bA\subset\dX$, we define
$(\bA-\bA):=\set{x-y}{x,y\in\bA}$.

\Proposition{\label{still.life.thm}}
{
 Let $\gk=\lam[\dK]^{-1}\cdot\chr{\dK}$, where $\dK\subset\dX$ is compact. 
Let $\bA\subset\dX$.  Suppose 
$(\lam[\bA]/\lam[\dK])\in\CO{s_0,b_0}$ and that $(\bA-\bA)
\subset\dK$.   Then $\chr{\bA}$ is a still life.
}
\bthmprf Let $\ba:=\chr{\bA}$.  
For any $x\in\bA$,  let $\dK_x := \set{x-k}{k\in\dK}$. 

\Claim{ For any $x\in\bA$,  \ $\bA\subset\dK_x$.}
\bclaimprf 
If $y\in\bA$, then $k:=x-y \in (\bA-\bA) \subset \dK$, so
$k\in\dK$, so $y = x-k\in\dK_x$.
\eclaimprf

  Let $\gs\circ(\gk*\ba)=\chr{\bS}$ and
 $\gb\circ(\gk*\ba)=\chr{\bB}$, for some $\bS,\bB\subset\dX$.

\Claim{$\bA\subset\bS$.}
\bclaimprf 
For any $x\in\bA$, \
$\gk*\ba(x) \ = \ \lam[\bA\intsct\dK_x]/\lam[\dK]
\ \eeequals{(*)} \ \lam[\bA]/\lam[\dK] \ \in \CO{s_0,b_0}\subset\CC{s_0,s_1}$,
where $(*)$ is by Claim 1.  Thus, $\gs\circ(\gk*\ba)(x) \ = \ 1$,
so $x\in\bS$.
\eclaimprf

\Claim{$\bB=\emptyset$.
 (Hence, it is vacuously true that $\bB\subseteq\bA$).}
\bclaimprf
For any $x\in\dX$, \ 
$\gk*\ba(x) \ = \ \lam[\bA\intsct\dK_x]/\lam[\dK]
\ \leq \lam[\bA]/\lam[\dK] \ < \ b_0$.

Thus, $\gb\circ(\gk*\ba)(x)=0$.\eclaimprf

Claims 2 and 3 satisfy the conditions of Prop.\ref{fixed.point.cond},
so $\ba$ is a still life.
\ethmprf

The main examples of Proposition \ref{still.life.thm}
are  balls with respect to
some norm on $\dX=\Real^D$.

\Proposition{\label{norm.ball.still.life}}
{
 Let $\norm{\bullet}{*}$ be a norm on $\dX$, and for any $r>0$, let
$\bigodot(r):=\set{x\in\dX}{\norm{x}{*}\leq r}$.
Let $\dK:=\bigodot(1)$, and let $\gk := \lam[\dK]^{-1}\cdot \chr{\dK}$.
 Suppose $s_0 \leq \frac{1}{2^D}$, and let
  $R := \min\{\sqrt[D]{b_0}, \ \frac{1}{2}\}$.

If $r<R$, \ 
$\bA \subseteq \bigodot(r)$, and $s_0\cdot\lam[\dK] \leq \lam[\bA]$
 then $\ba=\chr{\bA}$ is a still life.

In particular,
if $\bigodot(\sqrt[D]{s_0})\subseteq\bA\subseteq\bigodot(r)\subset\bigodot(R)$,
then $\ba=\chr{\bA}$ is a still life.
}

  Before proving Proposition \ref{norm.ball.still.life}, 
we give some examples. 
Let $\dX=\Real^2$ and suppose $s_0\leq \frac{1}{4} \leq b_0$,
so that $R=\frac{1}{2}$.
\bdesc
  \item[$\ell^1$ norm:] For any $r>0$, let 
$\bD(r)  :=  \set{x=(x_1,x_2)\in\dX} {|x_1|+|x_2| \leq r}$ be
the diamond of diameter $2r$. Let  $\gk := \frac{1}{2} \chr{\bD(1)}$. 
Then $\chr{\bD(r)}$ is a still life for any $r\in\CO{\sqrt{s_0},\frac{1}{2}}$.
 
 \item[$\ell^2$ norm:] For any $r>0$, let 
$\bB(r) := \set{x\in\Real^2}{|x|\leq r}$ be the disk of radius $r$.
 Let $\gk := \frac{1}{\pi} \chr{\bB(1)}$. 
Then $\chr{\bB(r)}$ is a still life for any $r\in\CO{\sqrt{s_0},\frac{1}{2}}$,

  \item[$\ell^\oo$ norm:] For any $r>0$, let 
$\bC(r)  :=  \CC{-r,r}^2$ be the square of sidelength $2r$.
Let  $\gk := \frac{1}{4} \chr{\bC(1)}$. 
Then $\chr{\bC(r)}$ is a still life for any $r\in\CO{\sqrt{s_0},\frac{1}{2}}$.
\edesc

\bthmprf[Proof of Proposition \ref{norm.ball.still.life}:]
\quad We will verify the conditions of Proposition \ref{still.life.thm}.
To see that $(\bA-\bA)\subset\dK$,
suppose $x,y\in\bA$.  Then $x,y\in\bigodot(r)$, so
$\norm{x-y}{*}  \leq   \norm{x}{*} + \norm{y}{*}
 \leq  r + r < 2R \leq  1$.  Thus,
 $(x-y)\in\dK$.

Also, $(\lam[\bA]/\lam[\dK]) \in \CO{s_0,b_0}$, because
\[
s_0 \quad\leeeq{(*)}\quad
\frac{\lam[\bA]}{\lam[\dK]}
 \quad  \leeeq{(\dagger)} \quad  \frac{\lam[\bigodot(r)]}{\lam[\dK]}
\quad \eeequals{(\ddagger)} \quad r^D\cdot \frac{\lam[\dK]}{\lam[\dK]}
 \quad = \quad r^D \quad \lt{(\diamond)} \quad b_0.
\]
Here, $(*)$ is because  $s_0\cdot\lam[\dK]\leq \lam[\bA]$ by hypothesis.
\ $(\dagger)$ is because  $\bA\subseteq\bigodot(r)$.
\ $(\ddagger)$ is because $\bigodot(r) \ = \ r\cdot\dK$, and $\lam$ is the
$D$-dimensional Lebesgue measure.  \ $(\diamond)$ is because
$r<R\leq \sqrt[D]{b_0}$.

In particular, if $\bA\supseteq \bigodot(\sqrt[D]{s_0})$, then $\lam[\bA]
 \  \geq \  \lam[\bigodot(\sqrt[D]{s_0})]
\  \ \eeequals{(*)} \ \ s_0\cdot \lam[\dK]$, \
where $(*)$ is because $\bigodot(\sqrt[D]{s_0}) \ = \ \sqrt[D]{s_0}\cdot\dK$.
\ethmprf

{\em Remark:} Proposition \ref{norm.ball.still.life} doesn't apply
if $s_0=b_0$.  However, the proof
can be extended to the special case  $s_0=\frac{1}{2^D} = b_0$
(eg. $s_0=\frac{1}{4} = b_0$ when $D=2$).\hfill$\diamondsuit$

\breath

  Let $\gam:\Real\into\Real^2$ be a smooth path.  If $w>0$,  the
{\dfn ribbon} of width $w$ around $\gam$ is the set
\[
  \sR(\gam,w)\quad:=\quad \set{x\in\Real^2}{|x-\gam(t)|\leq w/2, \ 
 \mbox{for some $t\in\Real$}}.
\]
 We assume that $\gam$ is an arc-length
parameterization --ie. $|\dot\gam|\equiv 1$.  
The {\dfn curvature} of $\sR(\gam,w)$ is the maximal
value of $|\ddot{\gam}(t)|$ for $t\in\Real$.
In particular, a {\dfn flat ribbon} is one with  curvature
0 ---in this case, $\gam$ is an affine function (ie. $\gam(t) = tv + x$
for some $x,v\in\Real^2$, with $|v|=1$).
We'll construct still lifes shaped like slowly curving ribbons
through $\Real^2$.

\begin{figure}
\begin{tabular}{cc}
{
\psfrag{B}[][]{$\ss$}
\psfrag{S}[][]{$\undS_1$}
\includegraphics[width=17em,height=16em]{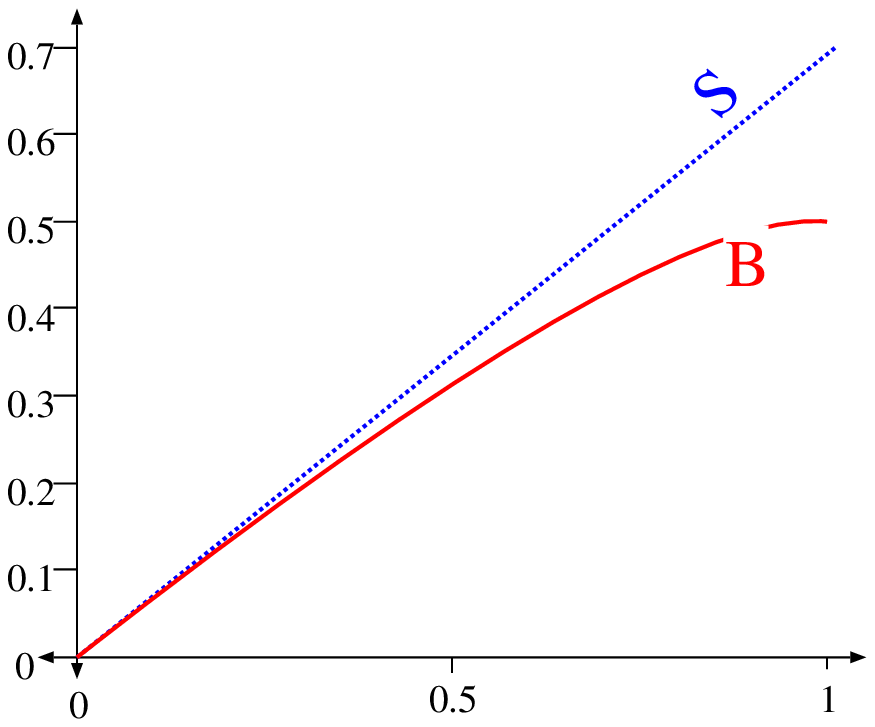}
}
 &
{\psfrag{B}[][]{\footnotesize $\dB_x$}
\psfrag{R}[][]{\footnotesize $\sR$}
\psfrag{x}[][]{\footnotesize $x$}
\psfrag{w}[][]{\footnotesize $w$}
\psfrag{w1}[][]{\footnotesize $\sqrt{1-w^2}$}
\psfrag{Area1}[][]{\footnotesize Area=$\frac{\pi}{4}-\arccos(w)/2$}
\psfrag{Area2}[][]{\footnotesize Area=$\frac{1}{2}w\cdot\sqrt{1-w^2}$}
\psfrag{theta}[][]{\footnotesize $\arccos(w)$}
\psfrag{theta2}[][]{\footnotesize $\frac{\pi}{2}-\arccos(w)$}
\includegraphics[scale=0.5]{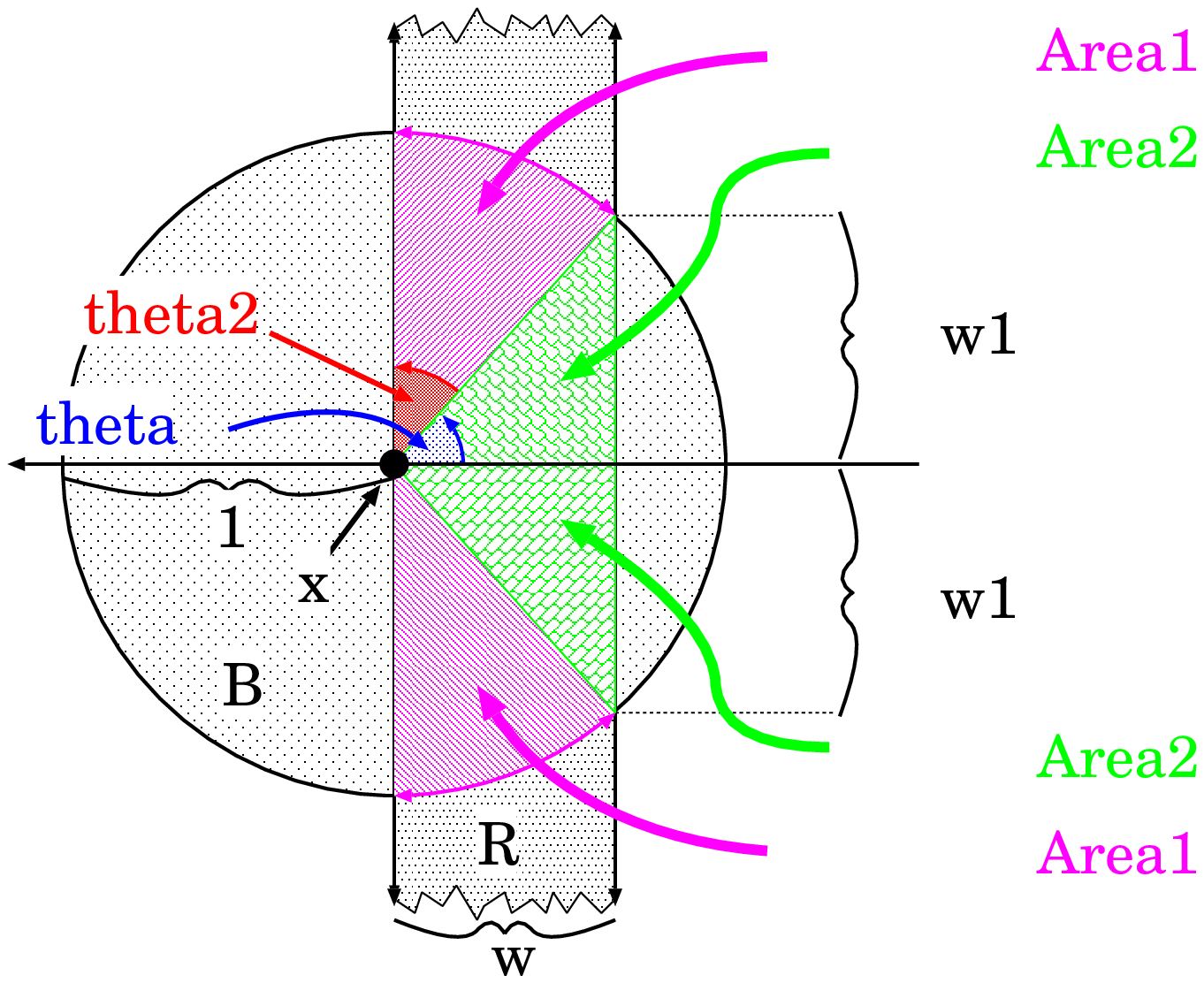}}\\
{\bf(A)} & {\bf(B)}
\end{tabular}
\caption{The functions $\undS_1$ and $\ss$ in
Proposition \ref{ribbon.still.life}.\label{fig:ribbon}}
\end{figure}

\Proposition{\label{ribbon.still.life}}
{
  Let $\dB:=\set{x\in\Real^2}{|x|\leq 1}$ be the
unit disk.  Let $\gk:=\frac{1}{\pi}\chr{\dB}$.
We define functions $\ss,\undS_1:\OC{0,1}\into\OC{0,1}$ by
\beq
\ss(w) &:=&
\frac{1}{2} + \frac{1}{\pi}
\lb( w\cdot\sqrt{1-w^2} - \arccos(w)\rb),\\
\And \ \undS_1(w) &:=&
1 + \frac{1}{\pi}
\lb( \frac{w}{2}\cdot\sqrt{4-w^2} - 2\arccos(w/2)\rb).
\qquad\mbox{\rm[Figure \ref{fig:ribbon}(A)]} 
\eeq
 Let $w\in\OC{0,1}$.
\bthmlist
  \item If $s_0  \leq  \ss(w)  \leq  b_0$
and $\undS_1(w)  \leq  s_1$, then any flat ribbon of width $w$ is a still life.

  \item   If $s_0  <  \ss(w)  <  b_0$
and $\undS_1(w)  <  s_1$, then any ribbon of width $w$
and small enough curvature is a still life.
\ethmlist
}
\bthmprf
  {\bf(a)} Let $\sR$ be a flat ribbon of radius $w$, and let
$\br:=\chr{\sR}$.  Suppose $x$ is a point on the boundary of $\sR$, and
$z$ is a point on the  centre line of $\sR$.

\Claim{\label{ribbon.still.life.claim}
{\bf(a)} \ $\gk*\br(x) = \ss(w)$. \qquad {\bf(b)} \  $\gk*\br(z) = \undS_1(w)$.}
\bclaimprf
Let $\dB_x=\set{b\in\Real^2}{|x-b|\leq 1}$.  Then
$\gk*\br(x) \  = \ \frac{1}{\pi} \lam[\dB_x\intsct\sR] = \ss(w)$, as shown in
Figure \ref{fig:ribbon}(B).  Likewise,
$\gk*\br(z) \  = \ \frac{1}{\pi} \lam[\dB_z\intsct\sR] = 2\cdot \ss(w/2)= \undS_1(w)$.
\eclaimprf
Now, $\gk*\br(y)$ is monotonically decreasing function of the
distance of $y$ from the centre line of $\sR$.  Thus,  
If $y\in\sR$ is any other point inside the ribbon, then we have 
\[
\hspace{-1em}  s_0 \ \  \leeeq{(h)} \ \  \ss(w) \ \  \eeequals{(a)} \ \  \gk*\br(x)
\ \   \leq \ \  \gk*\br(y) \ \  \leq \ \  \gk*\br(z)
\ \  \eeequals{(b)}\ \  \undS_1(w) \ \   \leeeq{(h)} \ \  s_1.
\]
Here, (h) is by hypothesis, \
(a) is by Claim \ref{ribbon.still.life.claim}(a),
and  (b) is by Claim \ref{ribbon.still.life.claim}(b).
Hence, $\gs[\gk*\br(w)]=1$.

 However, if $y\in\compl{\sR}$, then
$\gk*\br(y) \ < \ \gk*\br(x) \ \ \eeequals{(a)} \  \ \ss(w) \ \leeeq{(h)} \ b_0$,
where (a) is by Claim \ref{ribbon.still.life.claim}(a) and
 (h) is by hypothesis. Thus,  $\gb[\gk*\br(y)]=0$, as desired.

{\bf(b)}  In a slightly curving ribbon, 
Claim \ref{ribbon.still.life.claim}(a) 
will contain an error proportional to the curvature of $\sR$.
But if  $s_0 \ < \ \ss(w) \ < \ b_0$ and the curvature is small enough,
then we will still have $s_0\leq \gk*\br(x)$ for any $x\in\partial\sR$,
and $\gk*\br(y) < b_0$ for any $y\in\compl{\sR}$.

Likewise, Claim \ref{ribbon.still.life.claim}(b)
will contain an error proportional to the curvature of $\sR$.
But if $\undS_1(w) \ < \ s_1$  and the curvature is small enough,
we will still have $\gk*\br(z)\leq s_1$, and
thus, $\gk*\br(y)\leq s_1$  for any $y\in\sR$.
\ethmprf

{\em Remark:} 
 Figure \ref{fig:ribbon}(A) shows that we can find values of
$w\in\OC{0,1}$ satisfying $s_0 \ \leq \ \ss(w) \ \leq \ b_0$ for any
$s_0\leq b_0 \leq 0.5$.  Also, for most of this range, we have $\undS_1(w)
\approx \ss(w)$, so it is not hard to simultaneously achieve $\undS_1(w)\leq
s_1$, if $s_1$ is much larger that $b_0$.  For example, if $s_0\leq
0.25 \leq b_0$, and $0.3\leq s_1$, then a value of $w\approx 0.4$ will
suffice.\hfill$\diamondsuit$

\breath

  We can generalize Proposition \ref{ribbon.still.life} to $\Real^D$.
Let $\Life$ be a RealLife EA with convolution kernel $\gk\in\sK$.  We say
$\gk$ is {\dfn rotationally symmetric} if there is some function
$\kap:\CO{0,\oo}\into\CO{0,\oo}$ so that $\gk(x)=\kap|x|$ for all $x\in\dX$.
If $R>0$, let $\dB(R) := \set{x\in\Real^D}{|x|\leq R}$
and $\bb_R := \chr{\dB(R)}\in\lAX$.
If $r\in\CC{0,R}$, let $\dA(r,R) := \set{x\in\Real^D}{r\leq |x|\leq R}$
be the {\dfn bubble} with inner radius $r$ and outer radius $R$
(e.g. if $D=2$, then $\dA(r,R)$ is an annulus).
Let $\ba_{r,R} := \chr{\dA(r,R)}\in\lAX$.  

  Let $\sM\subset\dX$ be a smooth,
$(D-1)$-dimensional hypersurface, and let $w>0$.  The
{\dfn curtain} of width $w$ around $\sM$ is the set
\[
  \sC(\sM,w)\quad:=\quad \set{x\in\dX}{|x-m|\leq w/2, \ 
 \mbox{for some $m\in\sM$}}.
\]
The {\dfn curvature} of $\sC(\sM,w)$ is the maximal
curvature of any smooth path through $\sM$ obtained by
intersecting a 2-dimensional affine plane with $\sM$.
In particular, a {\dfn flat curtain} is one with  curvature
0 ---in this case, $\sM$ is a $(D-1)$-dimensional affine hyperplane. 
For example, if $r>0$, and $\sM$ is a sphere of radius $R_0 = r+w/2$,
then $\sC(\sM,w) \ = \ \dA(r,r+w)$ is a
bubble, whose curvature is inversely proportional to $r$.

\Proposition{\label{curtain.still.life}}
{
  Let $\Life$ have a rotationally symmetric  kernel $\gk$.
\bthmlist
  \item   For any $r\in\Real$, let $\dH_r :=\set{(x_1,\ldots,x_D)\in\dX}{x_1\leq r}$. 
  Define $h:\Real\into\Real$ by $h(r) := \int_{\dH_r} \gk(x) \ dx$,
and define $\ss,\undS_1:\OO{0,\oo}\into\OC{0,1}$ by $\ss(w):=h(w)-h(0)$ and
$\undS_1(w):= h(w/2)-h(-w/2)$.  Then $\ss$ and $\undS_1$ are
differentiably increasing, with the
following properties, for any $w>0$:
\bitem
  \item[{\bf[i]}] If $s_0  \leq  \ss(w)  \leq  b_0$
and $\undS_1(w)  \leq  s_1$, then any flat curtain of
 width $w$ is a still life.

  \item[{\bf[ii]}] If $s_0  <  \ss(w)  <  b_0$
and $\undS_1(w)  <  s_1$, then any curtain of width $w$
and small enough curvature is a still life.
\eitem

  \item Define $\undS_1:\CO{0,\oo}\into\CO{0,\oo}$ by
$\undS_1(R) := \int_{\dB(R)} \gk(x) \ dx$, for any $R\geq 0$.
There is a differentiably increasing function 
$\ss:\CO{0,\oo}\into\CC{0,1}$ so that,
for any $R>0$, if $s_0\leq \ss(R)< b_0$ and $\undS_1(R)\leq s_1$,
then $\bb_R$ is a still life.

\item Let $\Del:=\set{(r,R)\in\Real^2}{0<r<R}$.
There are differentiable functions 
$\ss,\undB_0,\barB_1,\undS_1:\Del\into\CC{0,1}$
(nondecreasing in $R$ and nonincreasing in $r$)
so that, for any $(r,R)\in\Del$,
if  $s_0  \leq \ss(r,R) < b_0$, \  $\undS_1(r,R)  \leq  s_1$,
and either  $\undB_0(r,R)  <  b_0$ or $b_1 < \barB_1(r,R)$
then $\ba_{r,R}$ is a still life.
\ethmlist
}
\bthmprf {\bf(a)}\quad Let $\sM$ be a hyperplane, let
$\sC=\sC(\sM,w)$ be a flat curtain of width $w$,
and let $\bc:=\chr{\sC}$.  If $x\in\partial\sC$, then
$\gk*\bc(x)=\ss(w)$.  If $z\in\sM$, then $\gk*\bc(z)=\undS_1(w)$.
These values are independent of the choice of points $x$ or $z$ (because
$\sM$ is translationally symmetric), and independent of
the orientation of $\sM$ (because $\gk$ is rotationally symmetric).
The proof of {\bf[i]} is now like Proposition \ref{ribbon.still.life}(a).
The proof of  {\bf[ii]} is like Proposition \ref{ribbon.still.life}(b).

 {\bf(b)}\quad Fix $R>0$, and let $x\in\partial\dB(R)$.
Define $\ss(R) := \gk*\bb_R(x)$.  This value does not depend on $x$,
because the function $\gk*\bb_R$ is rotationally symmetric, because
$\gk$ and $\bb_R$ are rotationally symmetric. 
Also observe that $\gk*\bb_R(0) = \int_{\dB(R)} \gk(x) \ dx
= \undS_1(R)$.

 The value of $\gk*\bb_R(y)$ is a nonincreasing function of $|y|$.
If $y\in\dB(R)$, then $|y|<|x|$, so
\[
\hspace{-2em}
  s_0 \ \  \leeeq{(h)} \ \  \ss(R) \ \  = \ \  \gk*\bb_R(x)
\ \   \leq \ \  \gk*\bb_R(y) \ \  \leq \ \  \gk*\bb_R(0)
\ \  = \ \  \undS_1(R) \ \  \leeeq{(h)} \ \  s_1,
\]
where (h) is by hypothesis.
Hence, $\gs[\gk*\bb_R(y)]=1$, as desired.
  If $y\in \compl{\dB(R)}$, then $|y|>|x|$, so
$\gk*\bb_R(y) \ \leq \ \gk*\bb_R(x) \  =   \ \ss(R) \ \lt{(h)} \ b_0$.
Thus,  $\gb[\gk*\bb_R(y)]=0$, as desired.

  {\bf(c)}   Fix $R>0$, and let $x\in\partial\dA(r,R)$ be a point with $|x|=R$.
Define $\ss(r,R) := \gk*\ba_{r,R}(x)$.  This value does not depend on $x$,
because the function $\gk*\ba_{r,R}$ is rotationally symmetric, because
$\gk$ and $\ba_{r,R}$ are rotationally symmetric.  We also define
\beq
\barB_1(r,R) \quad  :=\quad  \min_{x\in\dB(r)} \gk*\ba_{r,R}(x)
\quad  \leq \quad 
\undB_0(r,R) & :=&  \max_{x\in\dB(r)} \gk*\ba_{r,R}(x)  \\
 \And \quad 
\undS_1(r,R) &  :=&  \max_{x\in\dA(r,R)} \gk*\ba_{r,R}(x).
\eeq
If $y\in\dA(r,R)$, then \ $s_0 \  \leeeq{(h)} \  \ss(r,R) \  := \  \gk*\ba_{r,R}(x)
\   \leq \  \gk*\ba_{r,R}(y) \  \leq \  \undS_1(r,R)
 \   \leeeq{(h)} \  s_1$, so $\gs[\gk*\ba_{r,R}(y)]=1$, as desired.
  If $y\in\compl{\dA(r,R)}$, and $|y|>R$, then
$\gk*\ba_{r,R}(y) \ < \ \gk*\ba_{r,R}(x) \  =:   \ \ss(r,R) \ \lt{(h)} \ b_0$,
so  $\gb[\gk*\bb_r(y)]=0$, as desired.

If $y\in\compl{\dA(r,R)}$, and $|y|<r$, then either
$\gk*\ba_{r,R}(y) \ \leq \ \undB_0(r,R) < b_0$, or
$\gk*\ba_{r,R}(y) \ \geq \ \barB_1(r,R) > b_1$; either way,
$\gb[\gk*\ba_{r,R}(y)]=0$, as desired.
\ethmprf

{\em Remarks:} (a) Proposition \ref{curtain.still.life}(c) describes
two classes of bubble-shaped still lifes: those with `small' internal
cavity (ie. a small value of $r$), and those with large cavity (large $r$).
 The `small cavity' bubbles satisfy the condition $\barB_1(r,R)>b_1$,
because $r$ is small enough that $\gk*\ba_{r,R}(x) > b_1$ for all $x\in\dB(r)$.  The `large cavity' bubbles must instead
satisfy the condition $\undB_0(r,R)<b_0$.  Observe that $\ss(r,R)
< \undB_0(r,R)$ (by concavity); hence a large-cavity bubble requires
$s_0 \leq \ss(r,R) <  \undB_0(r,R) < b_0$ (which is impossible if $s_0=b_0$).

(b) Proposition \ref{curtain.still.life} can be extended to {\em Larger than Life} CA,
with two caveats. 

\qquad [i] {\em RealLife} EA with rotationally symmetric kernels are
isotropic, but LtL CA are inherently anisotropic due to lattice
effects.  Thus, the functions $\undB_0$, $\undS_1$, etc. will all have a
directional dependence, and some curtain directions will be `favoured' over
others.  

\qquad [ii]   The functions $\gk*\br$, $\gk*\bb_R$, and 
$\gk*\ba_{r,R}$ will decrease in discrete steps.  Thus, we must replace the
`boundary value' function $\ss$ with {\em two} functions, $\barS_0$
and $\undB_0$, measuring the value of $\gk*\br$ on the `inside edge'
and `outside edge' of the boundary, respectively.  The inequality
`$s_0 \leq \ss < b_0$' is then replaced two inequalities:
 `$s_0 \leq \barS_0$' and `$\undB_0 < b_0$'.  In general,
$\undB_0 < \barS_0$; hence, these two inequalities are simultaneously
satisfiable, even if $s_0=b_0$  (as is often true for the LtL CA studied in
\cite{Evans1,Evans2,Evans3,Evans4,Evans5}).\hfill$\diamondsuit$

\section{Robustness of still lifes in the Hausdorff metric
\label{S:hausdorff}}
 
If $\bX,\bY\subset\dX$ are closed sets, then the {\em Hausdorff
metric} from $\bX$ to $\bY$ is defined
\beqn
\label{haus.metric.defn}
  \HD{\bX,\bY}\quad:=\quad \frac{1}{2}\sup_{x\in\bX} \inf_{y\in\bY} \ d(x,y)
\ + \ 
\frac{1}{2} \sup_{y\in\bY} \inf_{x\in\bX} \ d(y,x).
\eeqn
We define the metric $d_*$ on $\lAX$ as follows:
for any $\ba,\bb\in\lAX$, if $\ba=\chr{\bA}$ and $\bb=\chr{\bB}$ (for some
$\bA,\bB\subset\dX$), then
  $d_*(\ba,\bb) \ := \ \norm{\ba-\bb}{1}
\ + \ \HD{\partial\bA,\partial\bB}$.

The main result of this section states conditions under which a still
life in $\lAX$ will be surrounded by a $d_*$-neighbourhood of other
still lifes.  First we need some machinery.  If $\ba\in\lAX$, we say
that $\ba$ is {\dfn $\gk$-smooth} if the function $\gk*\ba$ is $\sC^2$
on an open dense subset of $\dX$.  For example,
if $\gk\in\sC^2$, then {\em all} elements of $\AX$ are
$\gk$-smooth [because, if $\ba\in\AX$, then for any $c,d\in\CC{1...D}$,
 $\partial_c\, \partial_d \, (\gk*\ba)(x)=(\partial_c \,\partial_d \, \gk)*\ba(x)$ is defined and continuous at all $x\in\dX$].

 If $\gk$ is not smooth (eg.
 $\gk=\lam[\dK]^{-1}\chr{\dK}$ for  $\dK\subset\dX$),
then $\gk$-smoothness is still fairly common.
For example, call a subset
$\bA\subset\dX$ {\dfn smoothly open} if $\bA$ is open and $\partial\bA$
is a piecewise smooth manifold.  Thus, an open ball is smoothly open,
because its boundary is a sphere, and an open cube is smoothly open,
because its boundary is a union of $2D$ flat faces.

\Lemma{}
{
 Let $\dK$ be smoothly open and let $\gk:=\lam[\dK]^{-1}\chr{\dK}$.
If $\bA\subset\dX$ is also smoothly open,
 then $\ba:=\chr{\bA}$ is $\gk$-smooth.
}
\bthmprf (sketch) The boundaries of $\bA$ and $\dK$ can be decomposed
 into pieces which are the graphs of smooth functions from
 $\Real^{D-1}$ into $\Real$.  There is an open dense subset
 $\bY\subset\dX$ so that, if $y\in\bY$, then $\gk*\ba(y)$ can be
 expressed as a sum of the integrals under these graphs, over
 rectangular domains whose boundaries are determined by the
 coordinates of $y$.  Thus, these integrals are (locally) smooth
 functions of the coordinates of $y$. 
\ethmprf

 If $\ba$ is $\gk$-smooth, and $\gs\circ(\gk*\ba)=\chr{\bS}$ and $\gb\circ(\gk*\ba)=\chr{\bB}$ for some $\bS,\bB\subset\dX$, then let
\[
 M_\oo(\ba) \quad  := \quad
 \sup_{x\in\partial\bS \union\partial\bB} \ \frac{1}{|\grad(\gk*\ba)(x)|}.
\]
 
\begin{figure}
\centerline{
\includegraphics[angle=-90,scale=0.45]{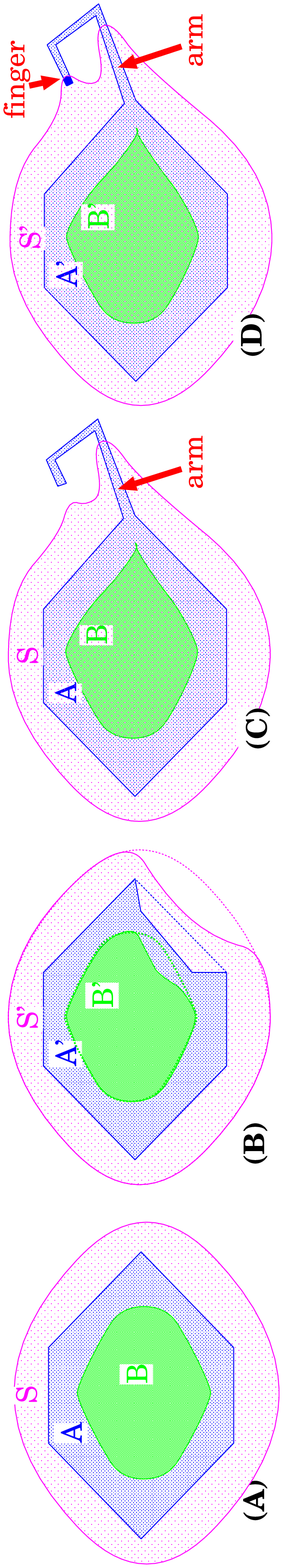}}
\caption{{\bf(A)} Theorem \ref{fixed.point.nhood} \
{\bf(B)} $\bB'\subset\bA'\subset\bS'$.
\ {\bf(C,D)} {\em RealLife} is not $d_*$-continuous.
\label{fig:boundary}}
\end{figure}

\Theorem{\label{fixed.point.nhood}}
{
 Suppose $\ba\in\lAX$ is $\gk$-smooth
still life, with  $M_\oo(\ba)<\oo$.
If $\cl{\bB}\subset \Int{\bA}$ and 
$\cl{\bA}\subset \Int{\bS}$ {\rm[Fig.\ref{fig:boundary}(A)]},
 then there is some $\eps>0$, so that
for any $\ba'\in\lAX$, if  $d_*(\ba,\ba')<\eps$
then $\ba'$ is also a still life.
}
\bthmprf
Our strategy is illustrated in Figure
\ref{fig:boundary}(B).  Suppose 
$\gb\circ(\gk*\ba')=\chr{\bB'}$.  If $\ba'$ is $\norm{\bullet}{1}$-close to
$\ba$, then  $\bB'$ will be $d_H$-close to $\bB$.  Thus, if 
$\cl{\bB}\subset \Int{\bA}$, and if $\bA'$ is  $d_H$-close to $\bA$,
then we'll have $\bB'\subset\bA'$.
Likewise, if $\gs\circ(\gk*\ba')=\chr{\bS'}$, then 
$\bS'$ will be $d_H$-close to $\bS$;  hence, if $\cl{\bA}\subset \Int{\bS}$,
then we'll have $\bA'\subset\bS'$, so the conditions
of Proposition \ref{fixed.point.cond} are satisfied.

  To realize this strategy, let $\alp:=\gk*\ba$, and
let $\bX:=\{x\in\dX$ \ ; \ $\alp$ differentiable and \ 
$\grad\alp(x)\neq 0\}$ (an open subset of $\dX$). 
Thus $\partial\bB\subset\bX$ and $\partial\bS\subset\bX$, because
$M_\oo<\oo$.

Now, $\alp\in\sC^2(\bX)$, so
define the $\sC^1$ vector field $\vV:\bX\into\Real^D$ by
 $\vV(x) \ := \ \frac{1}{|\grad\alp|^2}\grad\alp$,
and let $\sF:\bX\x\OO{-\tau,\tau}\into\bX$ be the flow induced by
$\vV$, which is well-defined in a
time-interval $\OO{-\tau,\tau}$ for some $\tau>0$.
In other words, for any $x\in\bX$, \
$\sF^0(x)=x$, and for any $t\in\OO{-\tau,\tau}$, \ 
$\partial_t \ \sF^t(x) \ = \ \vV\lb[\sF^t(x)\rb]$.

\Claim{\label{M.bound.C1}
For any $x\in\bX$, and $t\in\OO{-\tau,\tau}$, \ 
$\alp\lb[\sF^t(x)\rb] \ = \ \alp(x)+t$.}
\bclaimprf
  Let $\gam_x(t) := \alp\lb[\sF^t(x)\rb]$.  Thus, $\gam_x(0)=\alp(x)$.
If $t\in\OO{-\tau,\tau}$ and $\gam_x(t)=y$, then 
\[
\gam'_x(t) \quad = \quad \grad\alp(y)\bullet \vV(y)
\quad=\quad   \frac{\grad\alp(y)\bullet \grad\alp(y)}{|\grad\alp(y)|^2}
\quad=\quad 1.
\]
 Thus, $\gam_x(t) \ = \ \gam_x(0) + \int_0^t \gam_x'(r)  dr \ = \ 
\gam_x(0) + \int_0^t 1  dr \ = \  \gam_x(0)+t \ = \ \alp(x)+t$.
\eclaimprf

 Let $M:=M_\oo(\ba)$, and let $K:=\norm{\gk}{\oo}$.
Let $\del:=\norm{\ba-\ba'}{1}$.  Assume $\del<\tau/K$, so  $K\del<\tau$.

\Claim{\label{fixed.point.nhood.C2} {\bf(a)} \ \ 
$d_H(\partial\bS,\partial \bS') \ < \ M K \del + \sO\lb(K^2\del^2\rb)$.
\\ \mbox{}\qquad {\bf(b)} \ \ 
$d_H(\partial\bB,\partial \bB') \ < \ M K \del + \sO\lb(K^2\del^2\rb)$.}
\bclaimprf {\bf(a)}
$\bS=\alp^{-1}\CC{s_0,s_1}$, so
 $\partial\bS=\bC_0\disj\bC_1$, where
$\bC_0:=\alp^{-1}\{s_0\}$ and $\bC_1:=\alp^{-1}\{s_1\}$.
Likewise, if  $\alp'=\gk*\ba'$, then 
 $\partial\bS'=\bC'_0\disj\bC'_1$, where
$\bC'_i:=(\alp')^{-1}\{s_i\}$ for $i=1,2$.

Now $\norm{\alp-\alp'}{\oo}
\ \leq \ K\del$ by Lemma \ref{conv.lemma.2}(c).  Thus,
for any $x\in\dX$,  \ $|\alp(x)-\alp'(x)| < K\del$, so
$\bC'_0 \ = \ (\alp')^{-1}\{s_0\} \
 \subset \ \alp^{-1}\OO{s_0-K\del,s_0+K\del}$.
But $K\del<\tau$, so Claim \ref{M.bound.C1} implies that
\[
\alp^{-1}\OO{s_0-K\del,s_0+K\del} \quad = \quad 
\Union_{-K\del<t<K\del} \sF^t(\bC_0).
\]
  Thus, for any $c'\in\bC'_0$, there is some $c_0\in\bC_0$ and
some $t\in \OO{-K\del,K\del}$ such that 
\[
 c'\quad=\quad\sF^t(c_0)
\quad\eeequals{(\dagger)}\quad c_0 \ + \  t\vV(c_0) \ + \ \sO(t^2),
\]
where $(\dagger)$ is by Taylor's theorem. Thus,
\begin{eqnarray}
|c'-c_0| \ \ & = & \ \   \lb|t\vV(c_0) \ + \ \sO(t^2)\rb| 
\quad\leq\quad |t|\cdot\sup_{c\in\bC_0} |\vV(c)| \ + \  \sO(t^2)
\nonumber \\ &\leeeq{(*)}& \quad MK\del  \ + \  \sO(K^2\del^2).
\label{creeping.boundary.e1}
\end{eqnarray}
 $(*)$ is because $|t|<K\del$, and because 
$\D \sup_{c\in\bC_0} |\vV(c)| \ = \ 
\sup_{c\in\bC_0} \frac{1}{|\grad \alp (c)|} \ \leq \ M$.

Eqn.(\ref{creeping.boundary.e1}) means that
 $\D \inf_{c\in\bC_0} |c'-c| \ < \ MK\del  \ + \  \sO(K^2\del^2)$.
This is true for all $c'\in\bC'$; hence
$\D \sup_{c'\in\bC'_0} \inf_{c\in\bC_0} |c'-c| \ < \ MK\del  \ + \  \sO(K^2\del^2)$.
  Symmetric reasoning shows that
$\D \sup_{c\in\bC_0} \inf_{c'\in\bC'_0} |c-c'| \ < \ MK\del  \ + \  \sO(K^2\del^2)$. Hence, $d_H(\bC_0,\bC'_0) \ < \  MK\del  \ + \  \sO(K^2\del^2)$.

  By applying the same reasoning to $\bC_1:=\alp^{-1}\{s_1\}$ and
$\bC'_1:=(\alp')^{-1}\{s_1\}$, we can show that $d_H(\bC_1,\bC'_1) \ <
\ MK\del \ + \ \sO(K^2\del^2)$.  We conclude that
$d_H(\partial\bS,\partial \bS') \ < \ M K \del + \sO\lb(K^2\del^2\rb)$.
  The proof of {\bf(b)} is similar.
\eclaimprf
Now, let $\gam:=\min\lb\{d_H(\partial\bB,\partial\bA), \
d_H(\partial\bA,\partial\bS)\rb\}>0$.    Claim \ref{fixed.point.nhood.C2} 
yields some
$\del>0$ so that, if $\norm{\ba'-\ba}{1}<\del$, then
$d_H(\partial\bB',\partial\bB)<\gam/2$ and
$d_H(\partial\bS',\partial\bS)<\gam/2$.  Thus, if
$d_H(\partial\bA',\partial\bA)<\gam/2$, then
$\bB'\subseteq\bA'\subseteq\bS'$.  So let
$\eps:=\min\{\del,\gam/2\}$.  \ethmprf

\paragraph*{Remark:}
(a) $\Life$ is {\em not} $d_*$-continuous.  To see this, consider
Figure \ref{fig:boundary}(C), where $\bA$ is a hexagon with a long
`arm', and $\bS$ is an amorphous blob which contains the body of the
hexagon but not the arm.  Thus, $\Life(\chr{\bA})  =  \chr{\bA\intsct\bS}$
is the hexagon with most of the arm cut off.

  Now consider Figure \ref{fig:boundary}(D), where $\bA'$ is 
the same as $\bA$, but with a slightly longer arm, whose `finger'
rejoins the set $\bS'$.  Thus, $\Life(\chr{\bA'}) \ = \
\chr{\bA'\intsct\bS'}$ consists of the main hexagon, and also a
detached `finger' floating by itself.  The appearance of this
finger represents a discontinuous jump in the Hausdorff metric.
Thus, a $d_*$-continuous path in $\lAX$ from $\bA$ to $\bA'$ (by continuously
extending the arm) is mapped by $\Life$ into a $d_*$-discontinuous
path (where a finger suddenly appears).
Hence $\Life$ can't  be $d_*$-continuous along this path.

(b) If $b_0=s_0$, then parts of the boundaries of $\bB$,
$\bA$, and $\bS$ will generally coincide, so the hypotheses of Theorem
\ref{fixed.point.nhood} cannot be satisfied.
 
\subsection*{Conclusion}

This paper introduced the {\em
RealLife} family of Euclidean automata, and established their relationship to
{\em Larger than Life} CA.  Many questions remain.  For example, is the
converse to Theorem \ref{pointwise.convergence}(b) true?  That is: does 
a life form for {\em RealLife} imply the
existence of a  life form for LtL CA of sufficiently
large radius?  Also, \S\ref{S:rot.fixed.point} and \S\ref{S:hausdorff}
gave a variety of `existence' theorems for still lifes, but none for
other life forms.
Despite abundant empirical evidence, there are only a few rigorous
existence theorems regarding oscillators and bugs for LtL CA
\cite{Evans2,Evans3}, and as yet none for {\em RealLife}.

  Any compactly supported persistent structure (eg. an oscillator or bug) in a
cellular automaton must be eventually periodic, by the Pigeonhole Principle.
However, this is no longer true in Euclidean automata. Thus, {\em RealLife}
might possess {\em aperiodic} persistent structures, which are the
limits of a sequence of progressively longer-period oscillators or bugs
in progressively longer range LtL CA.  Is there an `evolution' theorem
for such structures, analogous to Theorem \ref{pointwise.convergence}(b)?

\ignore{
Is the analog of
Theorem \ref{branch.of.bugs} is true for {\em RealLife}, yielding
a parameterized family of bugs around a well-behaved oscillator?}

  Conway's {\em Life}  exhibits a complex
and subtle glider-based `physics', which makes possible
the construction of glider guns, glider reflectors and glider-based
logic, yielding machines capable of universal computation
and even self-replication \cite{BerlekampConwayGuy,RokaDurand}. 
Evans \cite{Evans5} 
has shown that at least one LtL CA ({\em Bosco's Rule}) exhibits similar 
universal computation (but not self-replication).
Does any {\em RealLife} EA contain  universal computers 
or self-replicators?

{\footnotesize
\bibliographystyle{alpha}
\bibliography{bibliography}
}

\end{document}